\documentstyle[amstex,amssymb,11pt]{article}

\begin{document}
    \pagestyle{plain} \setlength{\baselineskip}{1.3\baselineskip}
    \setlength{\parindent}{\parindent}

\title{{\bf Transverse measures, the modular class, and a cohomology pairing
for Lie algebroids}}  \author{Sam Evens and Jiang-Hua Lu\\ Department of
 Mathematics, University of Arizona, Tucson, AZ 85721
 USA\vspace{.2in}\\  Alan Weinstein\\Department of Mathematics,
 University of California, Berkeley, CA 94720 USA} \maketitle

 \newtheorem{thm}{Theorem}[section] \newtheorem{lem}[thm]{Lemma}
 \newtheorem{prop}[thm]{Proposition} \newtheorem{cor}[thm]{Corollary}
 \newtheorem{rem}[thm]{Remark} \newtheorem{exam}[thm]{Example}
 \newtheorem{nota}[thm]{Notation} \newtheorem{dfn}[thm]{Definition}
 \newtheorem{ques}[thm]{Question} \newtheorem{eq}{thm}

 \newcommand{\trace}{\mbox{Tr}} \newcommand{\del}{\partial}
 \newcommand{\rw}{\rightarrow} \newcommand{\lrw}{\longrightarrow}
 \newcommand{\rhu}{\rightharpoonup} \newcommand{\lhu}{\leftharpoonup}
 \newcommand{\Map}{\longmapsto} \newcommand{\qed}{\begin{flushright}
 {\bf Q.E.D.}\ \ \ \ \ \end{flushright} }
 \newcommand{\beqa}{\begin{eqnarray*}}
 \newcommand{\eeqa}{\end{eqnarray*}}

 \newcommand{\la}{\mbox{$\langle$}} \newcommand{\ra}{\mbox{$\rangle$}}
 \newcommand{\lala}{\mbox{$\la \! \la ~~~\ra \!\ra$}}

 \newcommand{\ot}{\mbox{$\otimes$}} \newcommand{\xa}{\mbox{$x_{(1)}$}}
 \newcommand{\xb}{\mbox{$x_{(2)}$}} \newcommand{\xc}{\mbox{$x_{(3)}$}}
 \newcommand{\ya}{\mbox{$y_{(1)}$}} \newcommand{\yb}{\mbox{$y_{(2)}$}}
 \newcommand{\yc}{\mbox{$y_{(3)}$}} \newcommand{\yd}{\mbox{$y_{(4)}$}}
 \renewcommand{\aa}{\mbox{$a_{(1)}$}}
 \newcommand{\ab}{\mbox{$a_{(2)}$}} \newcommand{\ac}{\mbox{$a_{(3)}$}}
 \newcommand{\ad}{\mbox{$a_{(4)}$}} \newcommand{\ba}{\mbox{$b_{(1)}$}}
 \newcommand{\bt}{\mbox{$b_{(2)}$}} \newcommand{\bc}{\mbox{$b_{(3)}$}}
 \newcommand{\ca}{\mbox{$c_{(1)}$}} \newcommand{\cb}{\mbox{$c_{(2)}$}}
 \newcommand{\cc}{\mbox{$c_{(3)}$}} \newcommand{\calH}{\mbox{$\cal
 H$}} \newcommand{\calS}{\mbox{$\cal S$}}
 \newcommand{\hyperH}{\mbox{$\Bbb H$}} \newcommand{\boldC}{\mbox{$\Bbb
 C$}} \newcommand{\reals}{\mbox{$\Bbb R$}}
 \newcommand{\zee}{\mbox{$\Bbb Z$}}

 \newcommand{\ts}{\mbox{$\sigma$}}
 \newcommand{\las}{\mbox{${}_{\sigma}\!A$}}
 \newcommand{\lasone}{\mbox{${}_{\sigma'}\!A$}}
 \newcommand{\ras}{\mbox{$A_{\sigma}$}}
 \newcommand{\rds}{\mbox{$\cdot_{\sigma}$}}
 \newcommand{\lds}{\mbox{${}_{\sigma}\!\cdot$}}

 \newcommand{\bb}{\mbox{$\bar{\beta}$}}
 \newcommand{\bg}{\mbox{$\bar{\gamma}$}}

 \newcommand{\id}{\mbox{${\em id}$}} \newcommand{\Fun}{\mbox{${\em
 Fun}$}} \newcommand{\End}{\mbox{${\em End}$}}
 \newcommand{\Hom}{\mbox{${\em Hom}$}}
 \newcommand{\ta}{\mbox{${\mbox{$\scriptscriptstyle A$}}$}}
 \newcommand{\ms}{\mbox{${\mbox{$\scriptscriptstyle M$}}$}}
 \newcommand{\ap}{\mbox{$A_{\mbox{$\scriptscriptstyle P$}}$}}
 \newcommand{\tx}{\mbox{$\mbox{$\scriptscriptstyle X$}$}}
 \newcommand{\pp}{\mbox{$\pi_{\mbox{$\scriptscriptstyle P$}}$}}
 \newcommand{\pg}{\mbox{$\pi_{\mbox{$\scriptscriptstyle G$}}$}}
 \newcommand{\asemi}{\mbox{$\ap \#_{\sigma} A^*$}}
 \newcommand{\dsemi}{\mbox{$A \#_{\Delta} A^*$}}

 \newcommand{\semi}{\mbox{$\times_{{\frac{1}{2}}}$}}
 \newcommand{\fk}{\mbox{${\frak k}$}} \newcommand{\fa}{\mbox{${\frak
 a}$}} \newcommand{\fd}{\mbox{${\frak d}$}}
 \newcommand{\ft}{\mbox{${\frak t}$}} \newcommand{\fg}{\mbox{${\frak
 g}$}} \newcommand{\fh}{\mbox{${\frak h}$}}
 \newcommand{\fn}{\mbox{${\frak n}$}} \newcommand{\fp}{\mbox{${\frak
 p}$}} \newcommand{\fbp}{\mbox{${\frak b}_{+}$}}
 \newcommand{\fbm}{\mbox{${\frak b}_{-}$}}
 \newcommand{\fnp}{\mbox{${\frak n}_{+}$}}
 \newcommand{\fnm}{\mbox{${\frak n}_{-}$}}
 \newcommand{\fgs}{\mbox{${\frak g}^*$}}
 \newcommand{\wg}{\mbox{$\wedge {\frak g}$}}
 \newcommand{\wgs}{\mbox{$\wedge {\frak g}^*$}}
 \newcommand{\wxl}{\mbox{$x_1 \wedge x_2 \wedge \cdots \wedge x_l$}}
 \newcommand{\wxk}{\mbox{$x_1 \wedge x_2 \wedge \cdots \wedge x_k$}}
 \newcommand{\wyl}{\mbox{$y_1 \wedge y_2 \wedge \cdots \wedge y_l$}}
 \newcommand{\wxkm}{\mbox{$x_1 \wedge x_2 \wedge \cdots \wedge
 x_{k-1}$}} \newcommand{\wxik}{\mbox{$\xi_1 \wedge \xi_2 \wedge \cdots
 \wedge \xi_k$}} \newcommand{\wxikm}{\mbox{$\xi_1 \wedge \cdots \wedge
 \xi_{k-1}$}} \newcommand{\wetal}{\mbox{$\eta_1 \wedge \eta_2 \wedge
 \cdots \wedge \eta_l$}}

 \newcommand{\winv}{\mbox{$(\wedge \fg_{1}^{\perp})^{\fg_1}$}}
 \newcommand{\wetak}{\mbox{$\eta_1 \wedge \cdots \wedge \eta_k$}}
 \newcommand{\gonep}{\mbox{$\fg_{1}^{\perp}$}}
 \newcommand{\wonep}{\mbox{$\wedge \fg_{1}^{\perp}$}}

 \newcommand{\db}{\mbox{$\fd = \fg \bowtie \fgs$}}
 \newcommand{\fds}{\mbox{${\scriptscriptstyle {\frak d}}$}}
 \newcommand{\fl}{\mbox{${\frak l}$}}

 \newcommand{\Gs}{\mbox{$G^*$}}
 \newcommand{\pis}{\mbox{$\pi_{\sigma}$}}
 \newcommand{\ea}{\mbox{$E_{\alpha}$}}
 \newcommand{\eb}{\mbox{$E_{-\alpha}$}} \newcommand{\Bm}{\mbox{$ {}^B
 \! M$}} \newcommand{\kBm}{\mbox{$ {}^B \! M^k$}}
 \newcommand{\Bb}{\mbox{$ {}^B \! b$}}
 \renewcommand{\epsilon}{\mbox{$\varepsilon$}}

 \newcommand{\cfg}{\mbox{$C(\fg \oplus \fgs)$}}
 \newcommand{\ps}{\mbox{$\pi^{\#}$}}
 \newcommand{\ppt}{\mbox{$\tilde{\pi}$}}
 \newcommand{\backl}{\mathbin{\vrule width1.5ex  height.4pt\vrule
 height1.5ex}}

 \newcommand{\bx}{\mbox{${\bar{x}}$}}
 \newcommand{\by}{\mbox{${\bar{y}}$}}
 \newcommand{\bz}{\mbox{${\bar{z}}$}}
 \newcommand{\pgs}{\mbox{${\pi_{\mbox{\tiny G}^{*}}}$}}

 \newcommand{\tp}{\mbox{$\varphi$}}
 \newcommand{\sn}{\mbox{$s_{\scriptscriptstyle N}$}}
 \newcommand{\tn}{\mbox{$t_{\scriptscriptstyle N}$}}
 \newcommand{\sm}{\mbox{$s_{\scriptscriptstyle M}$}}
 \newcommand{\tm}{\mbox{$t_{\scriptscriptstyle M}$}}
 \newcommand{\en}{\mbox{$\epsilon_{\scriptscriptstyle N}$}}
 \newcommand{\mem}{\mbox{$\epsilon_{\scriptscriptstyle M}$}}

 \newcommand{\qa}{\mbox{$Q_{A}$}}
 \newcommand{\pf}{\mbox{$\tilde{\pi}(df)$}}

 \newcommand{\ddp}{\mbox{$\delta^{'}_{\pi}$}}
 \newcommand{\pdp}{\mbox{$\partial_{\pi}$}}
 
 \begin{abstract}
 We show that every Lie algebroid $A$ over a manifold $P$ has a
 natural representation on the line bundle $ Q_A =
 \wedge^{top}A \otimes
 \wedge^{top} T^*P$. The line bundle $Q_A$ may be viewed as the 
Lie algebroid analog of the orientation bundle in topology, 
and sections of $Q_A$ may be
viewed as transverse measures to $A$. As a consequence, there is a
 well-defined class in
 the first Lie algebroid cohomology $H^1(A)$ called the modular class
 of the Lie algebroid $A$.  This is the same as the one
 introduced earlier by Weinstein using the Poisson structure on
 $A^*$. We show that there is a natural pairing between the Lie
 algebroid cohomology spaces  of $A$ with trivial coefficients and
 with coefficients in $Q_A$. This generalizes the pairing used in
 the Poincare duality of finite-dimensional Lie algebra cohomology.
 The case of holomorphic Lie algebroids is
 also discussed, where the existence of the modular class is connected
 with the Chern class of the line bundle $Q_A$.
 \end{abstract}

 \tableofcontents

 \section{Introduction}
 \label{sec_intro}

 The notion of Lie algebroids generalizes  that of both Lie algebras
 and tangent bundles.  Precisely, a {\bf Lie algebroid} over a smooth
 manifold $P$ is a vector bundle $A$ over $P$ together with 1) a Lie
 algebra structure on the space $\Gamma(A)$ of smooth sections of $A$,
 and  2) a bundle map $\rho: A \rightarrow TP$,  such that

 1) $\rho$ defines a Lie algebra homomorphism  from $\Gamma(A)$ to the
 space $\chi^1(P)$ of vector fields with the commutator Lie algebra
 structure, and

 2) for $f \in C^{\infty}(P)$ and $\omega_{1}, \omega_{2} \in
 \Gamma(A)$, the following derivation law holds:
 \[
 \{\omega_1, f \omega_2 \} = f \{\omega_1, \omega_2\} + (\rho
 (\omega_1) f) \omega_2.
 \]
 The map $\rho: A \rightarrow TP$ is  called the anchor map of the Lie
 algebroid $A$.  A {\bf representation} of a Lie algebroid $A$ over
 $P$ is a vector bundle $E$ over $P$, together with an
 $\reals$-bilinear map (we will deal with the complex case later)
 \[
 \Gamma(A) ~ \times ~ \Gamma(E) \lrw \Gamma(E): ~~ a ~ \ot~  s \Map
 D_a s,
 \]
 where $\Gamma(E)$ denotes the space of smooth sections of $E$, such
 that for any $a, b \in \Gamma(A), ~ s \in \Gamma(E)$ and $f \in
 C^{\infty}(P)$, \beqa & & (1) ~~ D_{fa}s ~ = ~ f D_a s;\\ & & (2)
 ~~ D_a(fs) ~ = ~ f D_a s + (\rho(a)f) s;\\ & &  (3) ~~ D_a(D_bs) -
 D_b(D_as) ~ = ~ D_{[a, b]} s.  \eeqa

 \bigskip
 Any Lie algebra is a Lie algebroid over a one point space, and its
 representations are representations of this Lie algebra.  The tangent
 bundle $TP$ is a Lie algebroid over $P$ with the identity map of $TP$
 as the anchor map. Representations  of $TP$ are vector bundles over
 $P$ with flat connections.  The trivial representation of a Lie
 algebroid $A$ is, by definition, the representation of $A$ on the
 trivial  line bundle over $P$ with the action given by
 \[
 D_a f ~ = ~ \rho(a) \cdot f, \hspace{.4in} a \in \Gamma(A), ~ f \in
 C^{\infty}(P).
 \] 

 \bigskip
 Our motivating example is, however, the cotangent  bundle Lie
 algebroid of a Poisson manifold: if $P$ is a Poisson manifold with
 the Poisson bivector  field $\pi$, the cotangent bundle $T^*P$ of $P$
 has a natural Lie algebroid structure, where the anchor map
 $\tilde{\pi}: T^*P \rightarrow TP$ is defined by
 \[
 \ppt(p): ~~T_{p}^{*}P \longrightarrow T_pP:~~ \alpha_p \Map  \alpha_p
 \backl \pi(p),
 \]
 and the Lie bracket of $1$-forms $\alpha$ and $\beta$ is given by
 \begin{eqnarray}
 \label{eq_bracket-on-one-forms}
 \{\alpha, \beta\} & = &d \pi (\alpha , \beta) ~ + ~ \ppt(\alpha)
 \backl d \beta ~ - ~ \ppt(\beta) \backl d \alpha\\
 \label{eq_bra-2}
 & = & - d \pi (\alpha, \beta) ~ + ~ L_{\tilde{\pi}(\alpha)} \beta ~ -
 ~ L_{\tilde{\pi}(\beta)} \alpha.
 \end{eqnarray}
 All the basic geometrical aspects of the Poisson structure $\pi$ such
 as its symplectic leaves and transversal Lie algebra structures to
 the symplectic leaves are reflected in the Lie algebroid  structure
 on $T^*P$. 

 Other examples of Lie algebroids are the gauge Lie algebroids of
 principal bundles, the boundary  Lie algebroids in $b$-calculus and
 the Weyl Lie algebroids in Fedosov quantization. 

 \bigskip
 Given a representation $E$ of a Lie algebroid $A$, one can define the
 Lie algebroid cohomology of $A$ with coefficients in $E$ (see Section
 \ref{sec_calculus}). When $A$ is a Lie algebra  or the tangent bundle
 $TP$, the Lie algebroid cohomology of $A$ is the Lie algebra
 cohomology or the de Rham cohomology of $P$. When $A = T^*P$ is the
 cotangent bundle Lie algebroid of a Poisson manifold $P$, the Lie
 algebroid cohomology of $A$ (with trivial coefficients)  is called
 the Poisson cohomology of $P$.

 \bigskip
 In this paper, we first construct, for an arbitrary Lie algebroid $A$
 over  $P$, an intrinsic representation of $A$ on the line bundle
 $Q_{A} = \wedge^{top}A \ot \wedge^{top}T^*P$. 
When $A$ is the subbundle of $TP$ tangent to a foliation ${\cal F}$,
sections of $Q_A$ are the transverse measures to the foliation, and 
we recover the Bott connection for this foliation.
Because of this, for a general Lie algebroid $A$, we may think of sections of
$Q_A$ as transverse measures to $A$.

\bigskip
We give two applications of the construction of the representation on
$Q_A$. 

 The first application concerns the modular class of $A$.  In
 \cite{we:modular}, the third author showed that there is a canonical class in
 the first Poisson cohomology space for each Poisson manifold $P$
 called the modular class of $P$, which measures the extent to which
 the Hamiltonian vector fields on $P$ fail to preserve volume forms on
 $P$ (when $P$ is orientable).  It is the semi-classical counterpart
 of the theory of modular automorphism groups for Von Neumann
 algebras. By using the Poisson structure on the dual bundle $A^*$ of
 a Lie algebroid $A$, he shows that there is also a canonical class,
 called the {\bf modular class} of $A$,  in the first Lie algebroid
 cohomology space of $A$ (with trivial  coefficients). We show that
 this modular class of $A$  can be directly constructed from the
 representation on $Q_A$. We also treat the case of holomorphic Lie
 algebroids in Section \ref{sec_holom}. We explain the obstruction
 for the existence of the modular class in this case in  terms of the
 Chern class of the line bundle $Q_A$. 

 When $A$ is the  cotangent bundle Lie algebroid $T^*P$ of a Poisson
 manifold $P$, we have $Q_A = (\wedge^{top}T^*P)^2$, and we show that
 there is, in fact, a representation of $A = T^*P$ on the "square
 root" $\wedge^{top}T^*P$ of $Q_A$. We also show that the cochain
 complex that calculates the Lie  algebroid cohomology of $T^*P$ with
 coefficients in $\wedge^{top}T^*P$ is isomorphic to chain complex on
 differential forms on $P$ introduced 
 by Koszul \cite{ko:crochet} and
 studied by Brylinski \cite{by:homo}.

 As the second application of the representation of $A$ on $Q_A$, we
 establish a pairing between the Lie algebroid cohomology of $A$ with
 trivial coefficients and that with coefficients in $Q_A$. This
 generalizes the pairing that gives Poincare duality for Lie algebra
 cohomology or for de Rham cohomology. As a special case, we get a
 pairing between Poisson cohomology and what we call the "twisted
 Poisson cohomology". This pairing may, however, be degenerate in
 general, as is seen by examples. The  problem of when it is
 non-degenerate is very interesting.  We hope to look at it in the
 future.

 \section{Differential calculus on Lie algebroids}
 \label{sec_calculus}

 In this section, we list some facts on the calculus on Lie algebroids
 that will be used in this paper. See  \cite {mac:lgla} for more
 details.

 \bigskip
 Let $A$ be a Lie algebroid over $P$ with anchor map $\rho$. For $k
 \geq  0$, let $\Gamma(\wedge^k A^*)$ be the  space of smooth sections
 of $\wedge^k A^*$. Define
 \[
 d_{\scriptscriptstyle A}: ~~ \Gamma(\wedge^{k-1}A^*) \lrw
 \Gamma(\wedge^k A^*):
 \]
 \begin{eqnarray}
 \label{eq_dA}
 (d_{\scriptscriptstyle A} \xi)  (a_1, ..., a_k) & = &  \sum_i
 (-1)^{i+1} \rho(a_i) (\xi (a_1, ..., \hat{a}_i, ..., a_k)) \\ &   + &
 \sum_{i < j} (-1)^{i+j} \xi(\{a_i, a_j\}, ..., \hat{a}_i, ...,
 \hat{a}_j, ..., a_k). \nonumber
 \end{eqnarray}
 It is well-defined and satisfies $d_{\scriptscriptstyle A}^2 =
 0$. The cohomology of $(\Gamma(\wedge^{\bullet} A^*),
 d_{\scriptscriptstyle A})$  is called the Lie algebroid cohomology of
 $A$ (with  trivial coefficients), and it is denoted by
 $H^{\bullet}(A)$.

 In the case when $A = TP$ is the tangent bundle of $P$, the Lie
 algebroid cohomology of $A$ is nothing but the de Rham cohomology of
 $P$.

 In the case when $A = T^*P$ is the cotangent bundle Lie algebroid of
 a Poisson manifold $(P, \pi)$, the Lie algebroid cohomology of $A$ is
 called the Poisson cohomology of $(P, \pi)$.

 \bigskip
 Suppose that $E$ is a representation of $A$. Let 
 \[
 \Gamma^k(A^*, E) ~ = ~ \Gamma(\wedge^k A^*) ~ \ot ~ \Gamma(E)
 \]
 be the space of ``$k$-forms" on $A$ with values in $E$. We  regard the
 action of $A$ on $E$ as defining a map 
 \[
 D: ~ \Gamma^{0}(A^*, E) = \Gamma(E) \lrw \Gamma^1(A^*, E).
 \]
 Using $d_{\scriptscriptstyle A}$, we can extend $D$ to a map 
 \[
 D: ~~ \Gamma^{k}(A^*, E) \lrw \Gamma^{k+1}(A^*, E)
 \]
 by the rule
 \[
 D(\xi ~ \ot ~ s) ~ = ~ d_{\scriptscriptstyle A} \xi ~ \ot ~ s  ~ + ~
 (-1)^{k} \xi ~ \ot ~ Ds,
 \]
 where $\xi \in \Gamma(\wedge^k A^*)$ and $s \in \Gamma(E)$.  Then
 $D^2 = 0$. The cohomology of $(\Gamma^{\bullet}(A^*, E), ~ D)$ is
 called the cohomology of $A$ with coefficients in $E$.  The operator
 $D$ satisfies
 \begin{equation}
 \label{eq_D-derivation}
 D(\xi \wedge \eta ~ \ot ~ s) ~ = ~ d_{\scriptscriptstyle A} \xi \wedge
 \eta ~ \ot ~  s ~ + ~ (-1)^{|\xi|} \xi \wedge D(\eta ~ \ot ~ s)
 \end{equation}
 for $\xi, \eta \in \Gamma(\wedge^{\bullet} A^*)$ and $s \in
 \Gamma(E)$.

 \bigskip
 The Lie bracket on the sections of $A$ can be extended to  the
 so-called Schouten bracket  $[ ~~ , ~~]$ on the space
 $\Gamma(\wedge^{\bullet} A) = \oplus_k \Gamma( \wedge^k A)$ of
 multi-sections of $A$. It is  characterized by the following
 properties: for $f, f_1, f_2 \in \Gamma(\wedge^0 A) = C^{\infty}(P),
 ~ a, a_1, a_2 \in \Gamma(A)$ and $X, Y, Z \in \Gamma(\wedge^{\bullet}
 A)$,
 \begin{eqnarray}
 \label{eq_sch-A-1}
 & & [X, Y] \in \Gamma(\wedge^{|X| + |Y| -1} A);\\
 \label{eq_sch-A-2}
 & & [a, ~ f] ~ = ~ \rho(a) \cdot f;\\
 \label{eq_sch-A-3}
 & & [a_1, ~ a_2] ~\text{is  the  Lie  bracket  in}~
 \Gamma(A);\\ 
 \label{eq_sch-A-4}
 & & [X, Y] =  - (-1)^{(|X|-1)(|Y|-1)} [Y, X];\\
 \label{eq_sch-A-5}
 & & [X, Y \wedge Z] = [X, Y] \wedge Z + (-1)^{(|X|-1)|Y|} Y \wedge
  [X, Z].
 \end{eqnarray}
 Explicitly, for $\xi \in \Gamma(\wedge^{|X| + |Y| -1} A^*)$, we have,
 \begin{equation}
 \label{eq_sch-A-explicit}
 (\xi, ~ [X, Y]) ~ = ~ (-1)^{(|X|-1)(|Y|-1)} i_{X}
 d_{\scriptscriptstyle A}i_{Y} \xi ~ - ~ i_Y d_{\scriptscriptstyle A}
 i_X \xi ~ + ~ (-1)^{|X|} i_{X \wedge Y} d_{\scriptscriptstyle A} \xi,
 \end{equation}
 where $i_X$ is the contraction operator by $X$, i.e.,
 \[
 (i_X \xi) (Y) ~ = ~ (\xi, ~ X \wedge Y).
 \]

 \bigskip
 We now recall the Lie derivative operators on $A$.  For a section  $a
 \in \Gamma (A)$, let 
 \[
 L_a: \Gamma(\wedge^k A)  \longrightarrow \Gamma(\wedge^k A), ~~~~  k
 \geq 0
 \]
 be given by
 \begin{equation}
 \label{eq_L1}
 L_a ( a_1 \wedge a_2 \wedge \cdots \wedge a_k)  ~ = ~ 
[a, ~ a_1 \wedge a_2 \wedge \cdots \wedge a_k]
~ = ~  \sum_{i=1}^k a_1
 \wedge \cdots \wedge [a, a_i]  \wedge \cdots \wedge a_k,
 \end{equation}
 where $a_i \in \Gamma(A), i = 1, 2, ..., k$. It is well-defined.  
We have, for $a, b \in \Gamma(A)$ and $X
 \in  \Gamma(\wedge^{\bullet} A)$,  
 \begin{eqnarray}
 \label{eq_la1}
 L_{[a, b]} & = & L_a L_b ~- ~ L_b L_a\\
 \label{la2}
 L_{fa}X & = & f L_a X ~ - ~ a \wedge (d_{\scriptscriptstyle A}f
 \backl X)\\
 \label{la3}
 L_a (fX) & = & f L_a X ~ + ~ (\rho(a) \cdot f)X.
 \end{eqnarray}
 We use the same letter $L_a$ to denote the operator on
 $\Gamma(\wedge^{\bullet} A^*)$ given by
 \begin{equation}
 \label{eq_L2}
 (L_a \xi, ~ X) ~ + ~ (\xi, ~ L_a X) ~ = ~ \rho(a) \cdot (\xi, ~ X),
 \end{equation}
 where $\xi \in \Gamma(\wedge A^*)$, and $X \in \Gamma(\wedge A)$.
 Then
 \begin{eqnarray}
 \label{eq_la4}
 L_a & = & d_{\scriptscriptstyle A} i_a ~ + ~  i_a
 d_{\scriptscriptstyle A}\\
 \label{eq_la5}
 L_{[a, b]} & = & L_a L_b ~- ~ L_b L_a\\
 \label{eq_la6}
 L_{fa}\xi & = & fL_a \xi ~ + ~ d_{\scriptscriptstyle A}f  \wedge
 i_a \xi\\
 \label{eq_la7}
 L_a(f \xi) & = & f L_a \xi ~ + ~ (\rho(a) \cdot f) \xi.
 \end{eqnarray}
 In particular, if $X \in \Gamma(\wedge^{top} A)$ and  $\xi \in
 \Gamma(\wedge^{top} A^*)$,  we have, for any $a \in \Gamma(A)$,
 \begin{eqnarray}
 \label{eq_la8}
 L_{fa}X & = & f L_a X ~ - ~ (\rho(a) \cdot f) X\\
 \label{eq_la9}
 L_a (fX) & = & f L_a X ~ + ~ (\rho(a) \cdot f) X\\
 \label{eq_la10}
 L_{fa}\xi & = & f L_a \xi ~ + ~ (\rho(a) \cdot f) \xi\\
 \label{eq_la11}
 L_{a} (f\xi) & = & f L_a \xi ~ + ~ (\rho(a) \cdot f) \xi.
 \end{eqnarray}

 \section{The representation of $A$ on $Q_A$ and the modular
 class of $A$}
 \label{sec_modular}

 Assume that a Lie algebroid $A$ acts on a real line bundle 
 $L$  over $P$. 
 We first assume that $L$ has a nowhere vanishing
 section; thus $L$ is trivial as a line bundle over $P$.

 Let $s$ be a nowhere vanishing  section of $L$. For a section $a$ of
 $A$, define $\theta_s(a) \in C^{\infty}(P)$ by
 \begin{equation}
 \label{eq_theta-s}
 D_a s ~ = ~ \theta_s(a) s.
 \end{equation}
 Clearly $\theta_s \in \Gamma(A^*)$. Moreover, it follows from 
 $D_{[a,b]} s ~ = ~ D_a D_b s ~ - ~ D_b D_a s$
 that 
 \[
 \theta_s([a, b]) ~ = ~ \rho(a) \cdot \theta_s(b) ~ - ~  \rho(b) \cdot
 \theta_s(a)
 \]
 for any two sections $a$ and $b$ of $A$. Thus $\theta_s$ is a
 $1$-cocycle with respect to $d_{\scriptscriptstyle A}$.  If $s_1$ is
 another nowhere vanishing section of $L$, and  if $s_1 = f_1 s$ for a
 nowhere vanishing function $f_1$,  we have
 \[
 \theta_{s_1} ~ = ~ \theta_s ~ + ~ d_{\scriptscriptstyle A}  (\log
 |f_1|) ~ \in  \Gamma (A^*).
 \]
  Thus $[\theta_{s_1}] = [\theta_s] \in H^1(A)$. We denote this class
 in $H^1(A)$ by $\theta_{L}$.

 Suppose that $L_1$ and $L_2$ are two line bundle  representations of
 $A$. Equip the tensor product line bundle $L_1  \ot L_2$ with the
 representation of $A$ given by 
 \[
 D_a (s_1 \ot s_2) ~ = ~ D_a (s_1) ~  \ot ~ s_2 ~ + ~ s_1 ~ \ot ~  D_a (s_2)
 \]
 for $a \in \Gamma(A)$. Then, assuming both $L_1$ and $L_2$ have
 nowhere vanishing sections, we have $\theta_{L_1 \otimes L_2} =
 \theta_{L_1} + \theta_{L_2}$. In particular, $ \theta_{L^2} = 2
 \theta_{L}$, where $L^2 = L \ot L$.

 \bigskip
 For a general line bundle $L$ that does not necessarily have a
 nowhere vanishing section, we define
 \[
 \theta_L ~ = ~ {\frac{1}{2}} \theta_{L^2} \in H^1(A).
 \]
 Since the square of any real line bundle over $P$ is trivial as a
 line bundle, the class $\theta_{L^2}$ is defined. By the remark
 above, $\theta_{L}$ is well-defined.

 \begin{dfn}
 \label{dfn_modular1}
 {\em The class $\theta_L \in H^1(A)$ is called the 
 {\bf characteristic class} of
 $A$ associated to the representation $L$.}
 \end{dfn}

 The following Proposition is immediate.

 \begin{prop}
 \label{prop_add}
 For any two line bundles $L_1$ and $L_2$ with $A$ representations,
 \[
 \theta_{L_1 \otimes L_2} ~ = ~ \theta_{L_1} ~ + ~ \theta_{L_2}.
 \]
 \end{prop}

Our notion of characteristic class is very much like that 
for flat line bundles.

 \begin{exam}
 \label{exam_trivial}
 {\em The characteristic class of $A$ associated to the trivial
 representation is zero.}
 \end{exam}

 \begin{exam}
 \label{exam_transformation}
 {\em Any Lie algebra homomorphism
 \[
 \rho: ~ \fg \lrw \chi^1(P)
 \]
 from a Lie algebra $\fg$ to the Lie algebra $\chi^1(P)$ of vector
 fields on $P$ defines a Lie algebroid structure on the trivial vector
 bundle $A = P \times \fg$ called the  {\bf transformation Lie
 algebroid}. The anchor map is $\rho$, now regarded as a bundle map
 from $P \times \fg$ to $TP$. The Lie bracket on the space $\Gamma(P
 \times \fg) \cong C^{\infty}(P, \fg)$ of smooth sections of $P \times
 \fg$ is given by
 \begin{equation}
 \label{eq_transformation-algebroid}
 \{ \bx, ~ \by\} ~ = ~ [\bx, ~ \by]_{\frak g} ~ + ~ \rho_{\bar{x}}
 \cdot \by ~ - ~ \rho_{\bar{y}} \cdot \bx,
 \end{equation}
 where the first term on the right hand side denotes the pointwise Lie
 bracket in ${\frak g}$, and the second term denotes the derivative of
 the ${\frak g}$-valued function $\by$ in the direction of the vector
 field $\rho_{\bar{x}}$.

 Suppose that $U$ is a vector space on which the Lie algebra $\fg$
 acts. Then  there is a representation of $A$ on the trivial vector
 bundle $P \times U$ given by
 \begin{equation}
 \label{eq_trans-on-U}
 D_{\bar{x}} \bar{u} ~ = ~  \rho_{\bar{x}} \cdot \bar{u} ~  + ~
 \bar{x} (\bar{u}),
 \end{equation}
 where $\bar{u} \in C^{\infty}(P, U)$ is a section of $P \times U$,
 $\rho_{\bar{x}} \cdot \bar{u}$ is the pointwise derivative of
 $\bar{u}$ in the direction of $\rho_{\bar{x}}$, and $\bar{x}
 (\bar{u})$ denotes the action of  $\bar{x} \in C^{\infty}(P, \fg)$ on
 $\bar{u} \in C^{\infty}(P, U)$ taken pointwise over $P$.  When $U$ is
 finite dimensional, we get a representation of  $A = P \times \fg$ on
 the trivial line bundle  $P \times (\wedge^{top}U)$. Let $\chi_U \in
 \fgs$ be the character of $\fg$ associated to $U$. Then the characteristic
 class of  the transformation Lie algebroid $A = P \times \fg$
 associated to the line bundle $P \times (\wedge^{top}U)$ is given by
 the constant section of $A^* = P \times \fgs$ defined by $\chi_U$.
 In particular, there is a representation of $A$ on the line bundle
 $\wedge^{top}A = P \times (\wedge^{top} \fg)$ coming from the adjoint
 action of $\fg$ on  $\wedge^{top} \fg$, and the corresponding characteristic
 class of $A$ is given by the constant section of $A^*$ defined by the
 adjoint character $\xi_0 \in \fgs$ of $\fg$.

 There is also a natural representation of the transformation Lie
 algebroid $A = P \times \fg$ on the tangent bundle $TP$ of $P$ via
 \begin{equation}
 \label{eq_trans-on-TP}
 D_{\bar{x}}V ~ = ~ [\rho_{\bar{x}}, ~ V] ~ + ~ \rho_{V \cdot
 {\bar{x}}},
 \end{equation}
 where $V$ is a vector field on $P$ and $V \cdot \bx$ denotes the Lie
 derivative of $\bx$ in the direction of $V$.  Correspondingly, there
 is a representation of $P \times \fg$ on the cotangent bundle $T^*P$
 satisfying
 \begin{equation}
 \label{eq_trans-on-T*P-1}
  \rho_{\bar{x}} (\alpha, ~ V) ~ = ~ (D_{\bar{x}} \alpha, ~ V) ~ + ~
 (\alpha, ~ D_{\bar{x}} V)
 \end{equation}
 for any $1$-form $\alpha$ on $P$.  Equivalently,
 \begin{equation}
 \label{eq_trans-on-T*P-2}
 (D_{\bar{x}} \alpha, ~ V) ~ = ~  \rho_{\bar{x}} (\alpha, ~ V) ~ - ~
 (\alpha, ~ [\rho_{\bar{x}}, ~ V] ~ + ~ \rho_{V \cdot {\bar{x}}})
 ~ = ~ (L_{\rho_{\bar{x}}} \alpha) (V) ~ - ~ (\alpha, ~ 
 \rho_{V \cdot \bar{x}}).
 \end{equation}
 Consequently, there is a representation of $A$ on the line bundle
 $\wedge^{top} T^*P$, as well as one on the density bundle
 $|\wedge^{top} T^*P|$ of $P$ (see \cite{bt:forms}).  Let $\mu$ be a non-vanishing density
 on $P$. Recall that the  divergence of a vector field $V$ on $P$ with
 respect to $\mu$ is defined to be the function $div_{\mu} V$ given by
 \[
 L_{V} \mu ~ = ~ (div_{\mu} V) \mu.
 \]
 The characteristic class of $A$ associated to  $\wedge^{top} T^*P$ or the
 density bundle  $|\wedge^{top}T^*P|$ is now given by the section of
 $A^*$ that associates to each constant section $x$ of $A$
 corresponding to $x \in \fg$  the function $div_{\mu} \rho_x$.  }
 \end{exam}

 \bigskip
 We now show that every Lie algebroid has an intrinsic
 line bundle representation. 

 \bigskip
 Let $A$ be a Lie algebroid over $P$ with anchor map $\rho$.  Consider
 the line bundle 
 \begin{equation}
 \label{eq_qa}
 Q_{\scriptscriptstyle A} ~ = ~  \wedge^{top} A ~ \ot ~  \wedge^{top}
 T^*P.
 \end{equation}
 For a section $a$ of $A$, define 
 \[
 D_a: ~~ \Gamma(\qa) \lrw \Gamma(\qa)
 \]
 by
 \begin{equation}
 \label{eq_d-on-qa}
 D_a(X \ot \mu) ~ = ~ L_a (X) ~ \ot ~ \mu ~ + ~ X ~ \ot ~
 L_{\rho(a)} \mu,
 \end{equation}
 where $X \in \Gamma(\wedge^{top} A), ~  \mu \in  \Gamma(\wedge^{top}
 T^*P)$, and $L_{\rho(a)}\mu$ denotes the Lie derivative of $\mu$ in
 the direction of $\rho(a)$.  It follows from (\ref{eq_la9}) and 
 (\ref{eq_la11}) that $D_a$ is well-defined. 

 \begin{thm}
 \label{thm_qa}
 The map 
 \[
 \Gamma(A) ~ \ot ~\Gamma(\qa) \lrw \Gamma(\qa): ~~
 a ~\ot ~ (X \ot \mu) \Map D_a (X \ot \mu)
 \]
 given by (\ref{eq_d-on-qa}) defines a representation of  $A$ on the
 line bundle $\qa$.
 \end{thm}

 \noindent
 {\bf Proof.} The proof is straightforward. The underlying idea is given in
 Remark \ref{rem_D-module}.  \qed

 \begin{dfn}
 \label{dfn_canonical-modular}
 {\em We define the {\bf modular class} of $A$ to be the
 characteristic class
 associated to the representation of $A$ on the  line bundle $\qa$,
 and we denote it by $\theta_{A}$ }
 \end{dfn}

 \begin{rem}
 \label{rem_D-module}
 {\em Alternatively, we can define $\qa$ as $Hom(\wedge^{top} A^*,
 \wedge^{top} T^*P)$. There is a Lie derivative action of $A$ on $
 \wedge^{top}(A^*)$, and a Lie derivative action of $A$ on $
 \wedge^{top}(T^*P)$ defined using the anchor map and the usual Lie
 derivative of vector fields.  Neither of these Lie derivative actions
 defines a representation of $A$, but we can combine these two Lie
 derivatives to get a representation of $A$ on $\qa$ by the formula
 \begin{equation}
 \label{eq_sam}
 (D_a \phi) (\xi) ~ = ~ L_{\rho(a)} (\phi(\xi)) ~ - ~  \phi(L_{a} \xi).
 \end{equation}
 When we calculate $(D_{fa} \phi)(\xi)$, each Lie derivative
 contributes a term $(\rho(a)f)\phi(\xi)$, and these two terms cancel:
 \beqa (D_{fa} \phi) (\xi) & = &  f\cdot L_{\rho(a)} (\phi(\xi))  ~ +
 ~ (\rho(a)f)\phi(\xi) ~ - ~ f\phi(L_{a} \xi) ~ - ~
 (\rho(a)f)\phi(\xi)\\ & = & f D_{a} \phi (\xi).  \eeqa Thus we obtain
 $D_{fa} = f D_a$.  In this interpretation, the representation of $A$
 on $\qa$ 
 is the difference between the Lie derivative action of
 $a$ on $\wedge^{top}A^*$ and the Lie derivative action on
 $\wedge^{top} T^*P$ defined via the anchor map.

 This construction can also be interpreted as follows: the  relations of
 the Lie derivative imply that we can give $\wedge^{top}A^*$ and
 $\wedge^{top}T^*P$ the structure of right representations of
 $A$. Explicitly, we set $r(a)\cdot s= ~-~L_a(s)$ and $r(f) \cdot s =
 fs$, where $f \in C^{\infty}(P)$ and $s$ is  a section of
 $\wedge^{top}A^*$ and use the analogous construction for
 $\wedge^{top}T^*P$.  A Lie algebroid has a functorially defined
 universal enveloping $D$-algebra, which is the associative algebra
 generated by functions on $P$ and sections of the Lie algebroid
 subject to the relations induced by the Lie algebroid relations (see
 \cite{bb:bb}).  When we say a module is a right representation of
 $A,$ we mean that it is a right module for the universal enveloping
 $D$-algebra.

 Let $D_P$ be the algebra of smooth differential operators on $P.$ It
 is the universal enveloping $D$-algebra for the tangent Lie algebroid
 $TP$. There is a standard construction in the theory of $D_P$-modules
 for changing a right $D_P$-module into a left one.  It proceeds by
 observing that the canonical bundle $\wedge^{top}T^*P$ is a right
 $D_P$-module and that for any right $D_P$-module $M$, there is a new
 left $D_P$-module  $Hom_{C^{\infty}(P)}(\wedge^{top}T^*P, M)$ with
 action given by the formula:
 \[
 ( a \cdot \phi) (\omega) ~ = ~ a (\phi(\omega)) -  \phi(L_a  \omega),
 \]
 where $a$ is a vector field on $P$ (see \cite{bo:dmod}). The action
 of $C^{\infty}(P)$ is by multiplication on $M$ and this makes
 $Hom_{C^{\infty}(P)}(\wedge^{top}T^*P, M)$ into a left
 $D_P$-module. Our construction of the representation of $A$ on $Q_A$
 is just the Lie algebroid analog of this construction applied to the
 right Lie algebroid module $\wedge^{top}T^*P$. It would be of
 interest to apply more serious ideas from $D$-module theory to Lie
 algebroids. We hope to  return to this subject in the future.  }
 \end{rem}

\begin{rem}
{\em
Another definition of the modular class begins with the observation of
Kontsevich \cite{ko:deformation} that a Lie algebroid $A$ can be viewed as a
supermanifold $A_s$, with $\Gamma(\wedge^{\bullet} A^*)$ as its algebra of
functions.  The differential $d_A$ defining Lie algebroid cohomology
can be considered as a vector field on $A_s$.  On the other hand, a
section $s_0$ of $Q_A$ determines in a natural way a volume element on
$A_s$.  The divergence of $d_A$ with respect to this volume
element is a function on $A_s$ which is in fact the section $\theta_{s_0}$
of $A^*$ representing the modular class.

}
\end{rem}

 \begin{rem}
 {\em When $P$ is not orientable, we can also replace
 $\wedge^{top}T^*P$ by the density bundle  $|\wedge^{top}T^*P|$, i.e.,
 define
 \[
 Q_A ~ = ~ \wedge^{top}A ~ \ot ~ |\wedge^{top}T^*P|.
 \]
 Then formula (\ref{eq_d-on-qa}) still defines a representation  of
 $A$ on $Q_A$.  }
 \end{rem}

 \begin{exam}
 \label{exam_tp}
 {\em For $A = TP$ and when $P$ is orientable, we have
 \[
 \qa ~ = ~ \wedge^{top} TP ~ \ot ~ \wedge^{top} T^*P ~ \cong ~
 Hom(\wedge^{top}T^*P, ~ \wedge^{top} T^*P).
 \]
 It has a canonical nowhere vanishing section given by the identity
 morphism. It is clear from (\ref{eq_sam}) that
 \[
 D_{V} s ~ = ~ 0
 \]
 for every vector field $V$ on $P$.  Thus the modular class of $TP$ is
 zero. When $P$ is not orientable, let ${\cal O}$ be the orientation
 bundle for $P$. Then $|\wedge^{top}T^*P| = \wedge^{top}T^*P \ot {\cal O}$,
 and thus $Q_A \cong {\cal O}$.
}
\end{exam}

\bigskip
 Because of Example \ref{exam_tp}, we can view $Q_A$ for
 a Lie algebroid $A$ as the {\bf orientation bundle of $A$}.

 \begin{exam}
 \label{exam_foliation}
 {\em When $A$ is an integrable subbundle of $TP$, the 
line bundle $Q_A$ is isomorphic to the top exterior power
of the conormal bundle to the corresponding foliation ${\cal F}$
(assuming the normal bundle is orientable - otherwise 
the twisting by the orientation bundle is needed), and thus
sections of $Q_A$ are transverse measures to ${\cal F}$.
The representation of $A$ on $Q_A$ is nothing but
the Bott connection of the foliation. It is thus appropriate
to call the modular class of $A$ the {\bf modular class of the 
foliation}.
}
\end{exam}

Because of Exam \ref{exam_foliation}, for a general Lie algebroid
$A$, we should think of sections of $Q_A$ as being ``{\bf transverse
measures to $A$}''.

\begin{exam}
 \label{exam_g}
 {\em When a Lie algebra $\fg$ is  considered as a Lie  algebroid over
 a one point space, its modular class  is given by the adjoint
 character $\xi_0 \in \fg^*$:
 \[
 (\xi_0, ~ x) ~ = ~ tr(ad_x), ~~~~ x \in \fg.
 \]
 }
 \end{exam}

\begin{exam}
 \label{exam_trans2}
 {\em Let $A = P \times \fg$ be the transformation Lie algebroid (see
 Example \ref{exam_transformation}). In this case, there is  a
 representation of $A$ on each of the line bundles $\wedge^{top}A$ and
 $\wedge^{top} T^*P$. The representation of  $A$ on $\qa$ is simply
 the tensor product of  these two. Let $\mu$ be a  nowhere vanishing
 density on $P$. Then the  modular class of $A$ is given by the
 section of $A^*$ that assigns to each section $\bar{x}$ of $A$ the
 function $(\xi_0, ~ \bar{x}) - div_{\mu}  (\rho_{\bar{x}})$, where
 $\xi_0 \in \fgs$ is the adjoint character of $\fg$.  }
 \end{exam}

 \section{The case of $T^*P$}
 \label{sec_t-star-P}

 In this section, we  treat the case of the cotangent bundle Lie
 algebroid $T^*P$ of a Poisson manifold $P$.  We will show that the
 modular class we define here is twice the modular class of $P$ as a
 Poisson manifold defined in \cite{we:modular}. We first
 recall  that definition.

 Consider the density bundle $|\wedge^{top}T^*P|$ of $P$. 
 Let $\mu$ be a nowhere  vanishing density
 on $P$. For a function $f$ on $P$, define $w_{\mu}(f) \in
 C^{\infty}(P)$ by
 \[
 L_{\tilde{\pi}(df)} \mu ~ = ~ w_{\mu}(f) ~ \mu.
 \]
 One then checks that $w_{\mu}(f)$ is a derivation with respect to
 $f$, so it is given by a vector field, which we denote by
 $w_{\mu}$. It also satisfies, for any $1$-forms $\alpha$ and $\beta$
 on $P$,
 \[
 w_{\mu}(\{\alpha, ~ \beta\}) ~ = ~  \ppt(\alpha) \cdot w_{\mu}(\beta)
 ~ - ~ \ppt(\beta) \cdot w_{\mu}(\alpha).
 \]
 (This is equivalent to $L_{w_{\mu}} \pi = 0$.) Thus, $w_{\mu}$
 defines a class in the first Poisson cohomology 
 $H^{1}_{\pi}(P)$ of $P$.  If $\mu_1 =
 f_1 \mu$ is another nowhere vanishing density on $P$, the two vector
 fields $w_{\mu}$ and $w_{\mu_1}$ are related by
 \[
 w_{\mu_1} ~ = ~ w_{\mu} ~ - ~ \ppt(d \log |f_1|).
 \]
 Thus $w_{\mu}$ and $w_{\mu_1}$ define the same  cohomology class in
 $H^{1}_{\pi}(P)$. 
 It is called  the modular class of $(P, \pi)$ in
 \cite{we:modular}.

 \bigskip
 We now look at the modular class $\theta_{T^*P}$  given by Definition
 \ref{dfn_canonical-modular}. The line bundle $\qa$ in this case is 
 \[
 \qa ~ = ~ (\wedge^{top}T^*P)^2.
 \]
 The representation of $T^*P$ on $(\wedge^{top}T^*P)^2$ is given by
 \[
 D_{\alpha} (\mu_1 \ot \mu_2) ~ = ~ \{\alpha, ~ \mu_1\} ~ \ot ~ \mu_2
 ~ + ~ \mu_1 ~ \ot ~ L_{\tilde{\pi}(\alpha)} \mu_2,
 \]
 where $\alpha$ is a $1$-form on $P, ~ \mu_1$ and $\mu_2$ are two top
 degree forms on $P$, and $\{\alpha, ~ \mu_1\}$ is
 the Schouten bracket between $\alpha$ and $\mu_1$. 

 We consider the special case when $\alpha = df$ is an exact $1$-form. 

 \begin{lem}
 \label{lem_df}
 For any differential $k$-form $\xi$ on $P$,
 \[
 \{df, ~ \xi\} ~ = ~ L_{\tilde{\pi}(df)} \xi.
 \]
 \end{lem}

 \noindent
 {\bf Proof.} Since the operators $L_{\tilde{\pi}(df)}$ and $L_{df}:
 \xi \mapsto \{df, ~ \xi\}$ have the same derivation properties with
 respect to the wedge product on  differential forms, it is enough to
 check the  cases when $k = 0$ and when $k = 1$. The case  $k = 0$ is
 obvious from the definition. Now let  $\xi$ be a $1$-form. From
 Formula (\ref{eq_bracket-on-one-forms}), we know that  for any
 $1$-form $\alpha$, \beqa \{\alpha, ~ \xi\} & = & (d
 i_{\tilde{\pi}(\alpha)}  ~ + ~ i_{\tilde{\pi}(\alpha)} d) \xi
 ~ - ~  \tilde{\pi}(\xi) \backl  d \alpha \\ & = &
 L_{\tilde{\pi}(\alpha)} \xi ~ - ~  \tilde{\pi}(\xi) \backl d \alpha.
 \eeqa In particular, when $d \alpha = 0$, we have
 \[
 \{\alpha, ~ \xi\} ~ = ~ L_{\tilde{\pi}(\alpha)} \xi.
 \]
 Thus $\{df, ~ \xi\} ~ = ~ L_{\tilde{\pi}(df)} \xi.$ \qed

 It now follows from the Lemma that 
 \[
 D_{df} (\mu_1 ~ \ot ~ \mu_2) ~ = ~ (L_{\tilde{\pi}(df)} \mu_1 )  ~
 \ot ~ \mu_2 ~ + ~ \mu_1 ~ \ot ~ (L_{\tilde{\pi}(df)} \mu_2)
 \]
 for any two top degree forms $\mu_1$ and $\mu_2$. 

 \bigskip
 The density bundle $S_P$ also satisfies $S_{P}^{2} =  (\wedge^{top}
 T^*P)^2$. Let $\mu$ be a nowhere vanishing section of $S_P$. Then
 $\mu \ot \mu = \mu^2$ is a nowhere vanishing section of
 $(\wedge^{top} T^*P)^2$.  For any $\alpha = df$, we have
 \[
 D_{df} ( \mu^2) ~ = ~ 2 (L_{\tilde{\pi}(df)} \mu ) ~ \ot ~ \mu.
 \]
 Since $L_{\tilde{\pi}(df)} \mu = w_{\mu}(f) \mu$, we have
 \[
 D_{df} ( \mu^2) ~ = ~ 2 w_{\mu}(f) \mu^2.
 \]
 Thus (recall the definition of $\theta_{\mu^2}$ from Section
 \ref{sec_modular}) we have
 \[
 \theta_{\mu^2} = 2 w_{\mu}.
 \]
 This shows that the modular class $\theta_{T^*P}$
 of the Lie algebroid $T^*P$ we define here is
 twice the modular class of $P$ as a Poisson manifold  defined in
 \cite{we:modular}.

 \begin{dfn} 
 \label{dfn_we}
 {\em The modular class $\theta_P$ of a Poisson manifold $P$ is
 defined to be $\theta_P = {\frac{1}{2}} \theta_{T^*P}$, i.e., one
 half of the modular class of the cotangent Lie algebroid $T^*P$.  }
 \end{dfn}

 Thus our definition for $\theta_P$ is the same as that in \cite{we:modular}.
 
 We now show that there is, in fact, a representation of  $T^*P$ on
 the line bundle $\wedge^{top}T^*P$.  We first state the following
 Proposition. The proof is straightforward, and we omit it.

 \begin{prop}
 \label{prop_sqrt}
 Assume that $M$ is a line bundle over $P$ such that a Lie algebroid
 $A$ with anchor map $\rho: A \rightarrow  TP$ acts on $M \ot M =
 M^2$. Then there is a representation of $A$ on $M$ as well that is
 defined as follows: let $a$ be a section of $A$, and let $s$ be a
 section of $M$.  Let $U$ be an open subset of $P$ over which $t$ is a
 nowhere vanishing section of $M$. Write $s = ft$ for some function
 $f$ defined on $U$. Define
 \begin{equation}
 \label{eq_D-sqrt}
 (D_a s)|_{U} ~ = ~ (\rho(a) \cdot f) t ~ + ~  {\frac{1}{2}}
 {\frac{D_a(t^2)}{t^2}} s.
 \end{equation}
 This is well-defined, and it defines a representation  of $A$ on $M$,
 whose square gives the original representation of $A$ on $M^2$.
 \end{prop}

 We now give the explicit representation of the  cotangent bundle Lie
 algebroid $T^*P$ on the line bundle $\wedge^{top}T^*P$, whose
 existence is  guaranteed by Proposition \ref{prop_sqrt}. Formula
 (\ref{eq_xu}) has also been independently discovered by Ping Xu
 in \cite{xu:homology}.

 \begin{thm}
 \label{thm_wedget*p}
 Let $P$ be a Poisson manifold with  the  Poisson bivector field
 $\pi$. For a $1$-form $\alpha$ and a top degree form $\mu$ on $P$,
 define 
 \begin{equation}
 \label{eq_d-alpha-mu}
 D_{\alpha} \mu ~ = ~ \{\alpha, ~ \mu\} ~ - ~ (\pi, ~ d \alpha) \mu,
 \end{equation}
 where $\{\alpha, ~ \mu\}$ denotes the Schouten bracket between
 $\alpha$ and $\mu$, and $(\pi, ~ d \alpha)$ denotes the pairing
 between $\pi$ and $d \alpha$.  Then (\ref{eq_d-alpha-mu}) defines a
 representation of $T^*P$ on $\wedge^{top}T^*P$. We also have
 \begin{eqnarray}
 \label{eq_another}
 D_{\alpha} \mu & = & L_{\tilde{\pi}(\alpha)} \mu ~ + ~  (\pi, ~ d
 \alpha) \mu\\
 \label{eq_xu}
 & = &  \alpha \wedge d i_{\pi} \mu.
 \end{eqnarray}
 \end{thm}

 \noindent
 {\bf Proof.} For any $f \in C^{\infty}(P)$, we have \beqa D_{f
 \alpha} \mu & = & \{f \alpha, ~ \mu\} ~ - ~ (\pi, ~  d(f \alpha))
 \mu\\ & = & f \{\alpha, ~ \mu\} ~ - ~ (\ppt(\alpha) \cdot f) \mu  ~ -
 ~ f (\pi, ~ d \alpha) \mu ~ -~  (\pi, ~ df \wedge \alpha) \mu\\ & = &
 f D_{\alpha} \mu ~ - ~ ((\pi, ~ \alpha \wedge df) +  (\pi, ~ df
 \wedge \alpha)) \mu\\ & = & f D_{\alpha} \mu;\\ D_{\alpha}(f \mu) & =
 & \{\alpha, ~ f \mu\} ~ - ~ (\pi, ~  d \alpha) f \mu\\ & = &
 (\ppt(\alpha) \cdot f) \mu ~ + ~ f D_{\alpha} \mu.  \eeqa If $\beta$
 is another $1$-form on $P$, we have 
 \beqa 
 D_{\alpha} D_{\beta} \mu ~
 - ~ D_{\beta} D_{\alpha} \mu & = & D_{\alpha}(\{\beta, ~ \mu\} -
 (\pi, ~ d\beta)\mu) ~ - ~ D_{\beta}(\{\alpha, ~ \mu\} - (\pi, d \alpha)
 \mu)\\ & = &\{\alpha, ~ \{\beta, ~ \mu\}\} ~ - ~ (\pi, ~ d \alpha) 
 \{\beta,~ 
 \mu\} \\ & & - ~ \ppt(\alpha) \cdot (\pi, ~ d \beta) \mu ~ -~ (\pi,~  d
 \beta) \{\alpha, ~ \mu\}  ~ + ~ (\pi, ~ d \beta) (\pi, ~ d \alpha) \mu\\ &
 & - ~ \{\beta, ~ \{\alpha, ~ \mu\}\} ~ + ~  (\pi,~ d \beta) \{\alpha,~ 
 \mu\} \\ & &  + ~ \ppt(\beta) \cdot (\pi, ~ d \alpha) \mu ~ + ~  (\pi, ~ d
 \alpha) \{\beta, ~ \mu\} ~ - ~  (\pi, ~ d\alpha) (\pi, ~ d \beta) \mu\\ &=&
 \{\{\alpha, ~ \beta\}, ~ \mu\} ~ - ~  \ppt(\alpha) \cdot (\pi, ~ d \beta)
 \mu ~ + ~ \ppt(\beta) \cdot (\pi, ~ d \alpha) \mu\\ & = & \{\{\alpha,
 ~ \beta\}, ~ \mu\} ~ - ~  (\pi, ~ d \{\alpha, ~ \beta\}) \mu\\ &=&
 D_{\{\alpha, ~ \beta\}} \mu.  \eeqa In the last step, we used the fact
 that
 \[
 \ppt(\alpha) \cdot (\pi, d \beta) ~ - ~ \ppt(\beta) \cdot (\pi, d
 \alpha) ~ = ~ (\pi, ~ d \{\alpha, \beta\}).
 \]
 Indeed,  \beqa {\text{ l.h.s}} & = & ([\ppt(\alpha), ~ \pi], ~ d \beta) ~ + ~
 (\pi, ~ L_{\tilde{\pi}(\alpha)} d \beta) ~ - ~  ([\ppt(\beta), ~ \pi],
 ~ d \alpha ) ~ - ~  (\pi, ~ L_{\tilde{\pi}(\beta)} d \alpha)\\ & = &
 -([\pi, ~ \ppt(\alpha)], ~ d \beta) ~ + ~  ([\pi, ~ \ppt(\beta)], ~ d
 \alpha) ~ + ~  (\pi, ~ d \{\alpha, \beta\})\\ & = & -(\ppt(d\alpha), ~
 d\beta) ~ + ~  (\ppt(d \beta), ~ d \alpha) ~ + ~ (\pi, ~  d \{\alpha,
 \beta\})\\ & = & (\pi, ~ d \{\alpha, \beta\})\\ & = & {\text {r.h.s}}.  \eeqa
 Here we used Formula (\ref{eq_bra-2}) and the fact that
 \[
 \ppt: ~~ \Gamma(\wedge^{\bullet} T^*P) \lrw \Gamma(\wedge^{\bullet}TP)
 \]
 intertwines the operators $d$ and $d_{\pi} = [\pi, ~ \bullet]$.

 To show that $D_{\alpha}\mu$ is also given by (\ref{eq_another}), we
 first observe that the right hand side of (\ref{eq_another}) is
 $C^{\infty}(P)$-linear in  $\alpha$. We also know from Lemma
 \ref{lem_df} that  (\ref{eq_another}) holds when $\alpha$ is exact.
 Thus it holds for an arbitrary $\alpha$. To show   (\ref{eq_xu}), we
 start with 
 \[
 \alpha \wedge i_{\pi} \mu ~ = ~ - i_{\tilde{\pi}(\alpha)} \mu
 \]
 and apply $d$ to both sides. We get
 \[
 d \alpha \wedge i_{\pi} \mu ~ - ~ \alpha \wedge d i_{\pi} \mu  ~ = ~
 - d i_{\tilde{\pi}(\alpha)} \mu,
 \]
 or,
 \[
 (\pi, ~ d \alpha) \mu ~ - ~ \alpha \wedge d i_{\pi} \mu  ~ = ~
 -L_{\tilde{\pi}(\alpha)} \mu.
 \]
 Formula (\ref{eq_xu}) now follows from (\ref{eq_another}).  \qed

 The characteristic class of the Lie algebroid $T^*P$ associated to  this
 representation on $\wedge^{top}T^*P$ is now exactly the same as the
 modular class of $P$ as a Poisson manifold defined in \cite{we:modular}.

 \bigskip
 We now consider the Lie algebroid cohomology  $H^{\bullet}(P, ~
 \wedge^{top}T^*P)$ of $T^*P$ with coefficients in $\wedge^{top}
 T^*P$. We will show that  it is isomorphic to the Poisson homology
 space  $H_{\bullet}(P)$ of Koszul  \cite{ko:crochet} and
 Brylinski  \cite{by:homo}.

 Recall (Section \ref{sec_calculus}) that the complex  that gives rise
 to the Lie algebroid cohomology  $H^{\bullet}(P, ~ \wedge^{top}T^*P)$
 is the space
 \[
 C_1 ~ = ~ \chi^{\bullet}(P) \ot \Gamma(\wedge^{top} T^*P)
 \]
 together with the operator
 \[
 \delta^{'}_{\pi}: ~~ \chi^k(P) \ot \Gamma(\wedge^{top} T^*P) \lrw
 \chi^{k+1}(P) \ot \Gamma(\wedge^{top} T^*P)
 \]
 given by
 \begin{equation}
 \label{eq_dp-1}
 \ddp: ~~ V \ot \mu \Map [\pi, ~ V] \ot \mu ~ + ~  (-1)^k V \wedge D
 \mu,
 \end{equation}
  where, for $\mu \in \Gamma(\wedge^{top} T^*P)$,
 \[
 D \mu \in \chi^1(P) \ot \Gamma(\wedge^{top} T^*P)
 \]
 is given by
 \[
 (D \mu ) (\alpha) ~ = ~ D_{\alpha} \mu  ~ = ~ \{ \alpha, ~ \mu\} ~ -
  ~ (\pi, d \alpha) \mu ~ = ~ L_{\tilde{\pi} \alpha} \mu ~ + ~ (\pi, d
  \alpha) \mu ~ = ~  \alpha \wedge d i_{\pi} \mu 
 \]
 for any $1$-form $\alpha$ on $P$ (see Theorem \ref{thm_wedget*p}).

 In \cite{ko:crochet}, Koszul (also see \cite{by:homo}) introduced
 the operator 
 \begin{equation}
 \label{eq_pp-1}
 \pdp ~ = ~ i_{\pi} d ~ - ~ d i_{\pi}: ~~ \Omega^{k+1}(P) \lrw
 \Omega^k (P).
 \end{equation}
 It satisfies $(\pdp)^2 = 0$. The homology of  $(\Omega^{\bullet}(P),
 \pdp)$ is called the  {\bf Poisson homology} of $(P, \pi)$ and it is
 denoted by $H_{\bullet}(P)$.

 Let $n = \dim P$, and for each $k = 0, 1, ..., n$, define
 \begin{equation}
 \label{eq_tau-p-1}
 \tau: ~~ \chi^k(P) \ot \Gamma(\wedge^n T^*P) \lrw  \Omega^{n-k}(P): ~~
  V \ot \mu \Map i_V \mu ~ = ~ V \backl \mu.
 \end{equation}
 It is clearly a vector space isomorphism from $C_1$ to
 $\Omega^{\bullet}(P)$.

 \begin{thm}
 \label{thm_iso-bry}
 For any $V \ot \mu \in C^{k}_{1} =  \chi^k (P) \ot \Gamma(\wedge^n
 T^*P)$,
 \begin{equation}
 \label{eq_tau-commute}
 \tau \ddp (V \ot \mu) ~ = ~ (-1)^{k+1} \pdp \tau  (V \ot \mu).
 \end{equation}
 Consequently,
 \[
 H^k (P, ~ \wedge^{top}T^*P) ~ \cong ~ H_{n-k} (P).
 \]
 \end{thm}

 \bigskip
 \noindent
 {\bf Proof.} Since all operators in (\ref{eq_tau-commute}) are local, 
 we only need
 to prove that it holds locally. Thus,  without loss of generality, we
 can assume that  $P$ is orientable. Let $\mu_0$ be a volume form of
 $P$.  We can then identify $\chi^{\bullet}(P)$ with $C_1$  by
 \[
 \chi^{\bullet}(P) \lrw C_1: ~~ V \Map V \ot \mu_0.
 \]
 The isomorphism $\tau: C_1 \rightarrow \Omega^{\bullet}(P)$ now
 becomes 
 \begin{equation}
 \label{eq_tau-p-0}
 \tau_{\mu_0}: ~ \chi^k (P) \lrw \Omega^{n-k} (P): ~  V \Map V \backl
 \mu_0 ~ = ~ i_V \mu_0.
 \end{equation}

 Let $\theta_0$ be the modular vector field associated  to
 $\wedge^{top}T^*P$ defined by $\mu_0$, i.e., for any $1$-form
 $\alpha$, 
 \begin{equation}
 \label{eq_theta-0}
 \theta_0(\alpha) \mu_0 ~ = ~ D_{\alpha} \mu_0  ~ = ~ L_{\tilde{\pi}
 (\alpha)} \mu_0 ~ + ~ (\pi, ~ d \alpha)  \mu_0
 \end{equation}
 (see Theorem \ref{thm_wedget*p}). Then the operator  $\ddp$ becomes 
 \begin{equation}
 \label{eq_dp-0}
 \delta^{'}_{\pi, \mu_0}: ~ \chi^k(P) \lrw \chi^{k+1} (P): ~  V \Map
 [\pi, ~ V] ~ + ~ \theta_0 \wedge V,
 \end{equation}
 and (\ref{eq_tau-commute}) becomes 
 \begin{equation}
 \label{eq_tau-commute-0}
 ([\pi, V] + \theta_0 \wedge V) \backl \mu_0 ~ = ~  (-1)^{k+1} (\pi
 \backl d(V \backl \mu_0) -  d((V \wedge \pi) \backl \mu_0))
 \end{equation}
 for $V \in \chi^k (P)$. 

 Introduce the operator 
 \begin{equation}
 \label{eq_b-0}
 b_{\mu_0}: ~ \chi^k(P) \lrw \chi^{k-1}(P): ~~ 
 (b_{\mu_0} V) \backl \mu_0 ~ = ~ (-1)^{k} d (V \backl \mu_0).
 \end{equation}
 Clearly $b_{\mu_0}^2 = 0$. 

 \begin{lem}
 \label{lem_b-0}
 For any $V_1, V_2 \in \chi^{\bullet}(P)$,
 \begin{equation}
 \label{eq_b-0-wedge}
 b_{\mu_0} (V_1 \wedge V_2) ~ = ~ b_{\mu_0} (V_1)  \wedge V_2 ~ + ~
 (-1)^{|V_1|} V_1 \wedge b_{\mu_0} (V_2) ~ + ~ (-1)^{|V_1|} [V_1, ~
 V_2];
 \end{equation}
 \end{lem}

 \bigskip
 \noindent
 {\bf proof.} 
 It follows from Formula (\ref{eq_sch-A-explicit}) for
 the Schouten bracket $[V_1, ~ V_2]$ that the two sides of
 (\ref{eq_b-0-wedge}) are equal when  paired with any for $\alpha$ of
 degree $|V_1| + |V_2| -1$.  \qed

 \begin{lem}
 \label{lem_b-0-pi}
 \begin{equation}
 \label{eq_b-0-pi}
 b_{\mu_0} \pi ~ = ~ \theta_0.
 \end{equation}
 \end{lem}

 \noindent
 {\bf Proof.} This is an immediate consequence of (\ref{eq_xu}).  \qed

 \noindent
  We now continue with the proof of Theorem \ref{thm_iso-bry}.
 We need to prove (\ref{eq_tau-commute-0}). Using again the operator 
  $b_{\mu_0}$ and its property (\ref{eq_b-0-wedge}), we
 see that 
 \beqa 
{\text {r.h.s. of}} ~(\ref{eq_tau-commute-0}) & = &
 (b_{\mu_0}(V \wedge \pi)  - b_{\mu_0} V \wedge \pi) \backl \mu_0\\ &
 = & ((-1)^k V \wedge b_{\mu_0} \pi + (-1)^k [V, \pi])  \backl \mu_0\\
 & = & (b_{\mu_0} \pi \wedge V ~ [\pi, V]) \backl \mu_0 \\ & = &
 (\theta_0 \wedge V + [\pi, V]) \backl \mu_0\\ & = & {\text {l.h.s.  of}} ~
 (\ref{eq_tau-commute-0}).  \eeqa This completes the proof of Theorem
 \ref{thm_iso-bry}.  \qed

 We now look at the case when $G$ is a Poisson Lie group.  Let $(G,
 \pi)$ be a connected Poisson Lie group with tangent Lie  bialgebra
 $(\fg, ~ \fgs)$ (see, for example,  \cite{sts:dressing} and
 \cite{lu-we:poi}).  Let $\xi_0 \in \fg$ and $x_0 \in \fg$ be
 respectively the  characters of the Lie algebras  $\fg$ and $\fgs$
 with respect  to their adjoint representations, i.e., for any $x \in
 \fg$ and  $\xi \in \fgs$, \beqa (\xi_0, ~ x ) & = & tr(ad_x)\\ (x_0,
 ~ \xi) & = & tr(ad_{\xi}).  \eeqa Denote respectively by $x_{0}^{l}$
 and $x_{0}^{r}$ the left and right invariant vector fields on $G$
 whose values at the identity element $e$ are $x_0$. Denote by
 $\xi_{0}^{r}$ the right invariant $1$-form on $G$ whose value at $e$
 is $\xi_0$. Recall that the  right dressing vector field 
 \cite{lu-we:poi} \cite{sts:dressing}
 defined by $\xi_0$ is the vector field
 \[
 \sigma_{\xi_0} ~ = ~ \ppt (\xi_{0}^{r}).
 \]

 \begin{prop}
 \label{prop_plg}
 Let $\mu$ be a right invariant volume form on $G$. The  vector field
 $\theta_{\mu}$ on $G$ defined by $\mu$:
 \[
 (\theta_{\mu}, ~ \alpha)  ~:= ~ {\frac{D_{\alpha}\mu}{\mu}}
 \]
 is given by
 \begin{equation}
 \label{eq_right}
 \theta_{\mu} ~ = ~ {\frac{1}{2}} (x_{0}^{l} ~ + ~  x_{0}^{r} ~ - ~
 \sigma_{\xi_0}).
 \end{equation}
 Similarly, if $\mu^{l}$ is a left invariant volume form on $G$.  The
 vector field $\theta_{\mu^l}$ on $G$ is given by
 \begin{equation}
 \label{eq_left}
 \theta_{\mu^l} ~ = ~ {\frac{1}{2}} (x_{0}^{l} ~ + ~ x_{0}^{r} ~ + ~
 \sigma_{\xi_0}).
 \end{equation}
 \end{prop}

 \noindent
 {\bf Proof.} Let $\alpha$ be a right invariant $1$-form on $G$.  The
 Schouten bracket $\{ \alpha, ~ \mu\}$ is also  right invariant, and
 is equal to $(x_{0}^{r}, \alpha) \mu$.  Consider the vector field $v$
 on $G$ defined by
 \[
 (v, ~ \alpha) ~ = ~ {\frac{1}{2}}(x_{0}^{r} - x_{0}^{l} ~ + ~
 \sigma_{\xi_0}, ~ \alpha) ~ - ~ (\pi, ~ d \alpha),
 \]
 where $\alpha$ is a right invariant $1$-form on $G$. It remains to
 show that $ v  =  0$.
 Since $\pi$ is multiplicative and since $\xi_0$ is $Ad_G$-invariant,
 we know that $v$ is multiplicative \cite{lu:thesis}, i.e.,
 for $g, h \in G$,
 \[
 v(gh) ~ = ~ l_g v(h) ~ + ~ r_h v(g).
 \]
 Since $v(e) = 0$, the fact that $v = 0$ follows from  the fact that
 the linearization of $v$ at $e$ is zero.  The case for $\mu^l$  is
 similarly proved.  \qed

 \begin{rem}
 \label{rem_modular}
 {\em 1) The vector fields $x_{0}^{l}$ and $x_{0}^{r}$ are Poisson
 vector fields because $x_0$ is a  character of $\fgs$.  Let $f_0$ be
 the function on $G$ defined by 
 \[
 f_0 (g) ~ := ~ \det(Ad_g ~ \mbox{on} ~ \wedge^{top} \fg).
 \]
 Then $d(\log f_0) = \xi_{0}^r$. Thus 
 \[
 \sigma_{\xi_0} ~ = ~ \ppt(d(\log f_0))
 \]
 is a Hamiltonian vector field.

 2) If we identify the cotangent bundle $T^*G$ with the trivial bundle
 $G \times \fgs$ by right translations, the Lie algebroid structure on
 $T^*G$ becomes that of the transformation Lie algebroid defined by
 the infinitesimal right dressing action of $\fgs$ on $G$ (see
 \cite{lu:thesis}). Proposition \ref{prop_plg} can then also be proved
 as a corollary of Example \ref{exam_trans2}.  }
 \end{rem}

 Assume now that $H \subset G$ is a connected and closed  Poisson Lie
 subgroup of $G$ \cite{sts:dressing}. Recall that  this is equivalent
 to ${\frak h}^{\perp} \subset \fgs$  being an ideal in $\fgs$, where
 \[
 {\frak h}^{\perp} ~ = ~ \{ \xi \in \fgs: ~ (\xi, x) = 0 ~~  \forall x
 \in \fh\}.
 \]
 The quotient space $G/H$ has a unique Poisson structure such that the
 projection map 
 \[
 j: ~ G \rightarrow G/H: ~~ g \Map gH
 \]
 is a Poisson map. Assume that there exists a $G$-invariant  volume
 form $\mu_0$ on $G/H$. We wish to determine the modular vector field
 $\theta_{\mu_0}$ on $G/H$ defined by $\mu_0$.

 Set $\tilde{\mu}_0 = j^* \mu_0$. It is a left invariant $l$-form on
 $G$, where $l = \dim (G/H)$. Choose any $\mu_1 \in \wedge^{n-l} \fgs$,
 where $n = \dim G$, such that $\mu =  \mu_1 \wedge \tilde{\mu}_0$ is a
 left invariant volume form on $G$. We have shown that the modular
 vector field $\theta_{\mu}$ on $G$ defined by $\mu$ is given by
 (\ref{eq_left}). 

 \begin{prop}
 \label{prop_quotient}
 The modular vector field $\theta_{\mu_0}$ on $G/H$ defined by $\mu_0$
 is given by
 \begin{equation}
 \label{eq_modular-quotient}
 \theta_{\mu_0} ~ = ~ j_{*} (\theta_{\mu}) ~ = ~  {\frac{1}{2}}  j_*
 (x_{0}^{l} ~ + ~ x_{0}^{r} ~ + ~  \sigma_{\xi_0}),
 \end{equation}
 i.e., it is the projection to $G/H$ by $j$ of the  modular vector
 field $\theta_{\mu}$ on $G$ defined by  $\mu = \mu_1 \wedge
 \tilde{\mu}_0$.
 \end{prop}

 \noindent
 {\bf Proof.} Let $\gamma_1, ..., \gamma_l$ be a basis for  ${\frak
 h}^{\perp} \subset \fg^*$. Denote by the same letters the
 corresponding left invariant $1$-forms on $G$. Assume that $\mu_1 =
 \xi_1 \wedge \xi_2 \wedge \cdots \wedge \xi_t$, where $t = \dim H$
 and the $\xi_i$'s are in $\fgs$ and denote also the corresponding
 left invariant $1$-forms on $G$. 

 Let $\alpha$ be an arbitrary $1$-form on $G/H$. Let $\tilde{\alpha} =
 j^* \alpha$ be the pull-back of $\alpha$ to $G$ by the projection map
 $j$. Then
 \[
 \{ \tilde{\alpha}, ~ \mu \} ~ = ~ \{\tilde{\alpha}, ~ \mu_1 \} \wedge
 \tilde{\mu}_0 ~ + ~ \mu_1 \wedge \{\tilde{\alpha}, ~ \tilde{\mu}_0\}.
 \]
 Now,
 \[
 \{\tilde{\alpha}, ~ \mu_1\} ~ = ~ \sum_{j=1}^{t}  \xi_1 \wedge \cdots
 \wedge \xi_{j-1} \wedge \{\tilde{\alpha}, ~ \xi_j \}  \wedge
 \xi_{j+1} \wedge \cdots \wedge \xi_t.
 \]
 Write $\tilde{\alpha} ~ = ~ \sum_{i=1}^{l} f_i \gamma_i$, where the
 $f_i$'s are functions on $G$. Then for each $j = 1, ..., t$,
 \[
 \{\tilde{\alpha}, ~ \xi_j\} ~ = ~ \sum_i(f_i \{\gamma_i, ~ \xi_j\} ~
 - ~ (\tilde{\pi}(\xi_j) \cdot f_i) \gamma_i).
 \]
 Since ${\frak h}^{\perp}$ is an ideal of $\fgs$, we have
 $\{\gamma_i, ~ \xi_j\} \in {\frak h}^{\perp}$. Thus
 \[
 \{\tilde{\alpha}, ~ \xi_j\} \wedge \tilde{\mu}_0 ~ = ~ 0, ~~~\forall
 j,
 \]
 and
 \[
 \{\tilde{\alpha}, ~ \mu_1\} \wedge \tilde{\mu}_0 ~ = ~ 0.
 \]
 Hence
 \[
 \{\tilde{\alpha}, ~ \mu\} ~ = ~ \mu_1 \wedge \{\tilde{\alpha}, ~
 \tilde{\mu}_0 \} ~ = ~  \mu_1 \wedge j^* \{\alpha, ~ \mu_0 \}.
 \]
 Therefore, using $\pi_G$ and $\pi$ to denote respectively  the
 Poisson bi-vector fields on $G$ and on $G/H$,  we have from the
 definition that \beqa \theta_\mu (\tilde{\alpha}) \mu & = &
 \{\tilde{\alpha}, ~ \mu_1 \wedge \tilde{\mu}_0 \} ~ - ~  (\pi_{G}, ~
 d \tilde{\alpha}) \mu_1 \wedge \tilde{\mu}_0\\ & = & \mu_1 \wedge j^*
 \{\alpha, ~ \mu_0 \} ~ - ~ j^* (\pi, ~ d \alpha) \mu_1 \wedge
 \tilde{\mu}_0\\ & = & \mu_1 \wedge j^*(\theta_{\mu_0} (\alpha) \mu_0)
 \\ & = & j^*(\theta_{\mu_0} (\alpha)) \mu.  \eeqa Thus
 \[
 \theta_\mu (\tilde{\alpha}) ~ = ~ j^*(\theta_{\mu_0} (\alpha)).
 \]
 This shows that $j_{*} \theta_{\mu}$ is well-defined and that
 \[
 \theta_{\mu_0} ~ = ~ j_{*} \theta_{\mu} ~ = ~  {\frac{1}{2}}  j_*
 (x_{0}^{l} ~ + ~ x_{0}^{r} ~ + ~ \sigma_{\xi_0}).
 \]
 \qed

 \begin{exam}
 \label{exam_flag}
 {\em Let $K$ be a compact semi-simple Lie group and let $T$ be a
 maximal torus of $K$. There is a naturally defined Poisson structure
 on $K$ making it into a Poisson Lie group \cite{lu-we:poi}. This
 Poisson structure vanishes at points in $T$. Thus it descends to a
 Poisson structure, called the Bruhat-Poisson structure, on the  flag
 manifold $K/T$. It follows from Proposition \ref{prop_quotient} that
 the modular vector field of the Bruhat-Poisson structure defined by a
 $K$-invariant volume form on $K/T$ is the vector field defined by
 $2iH_{\rho}$ in the Lie algebra ${\frak t}$ of $T$, where $\rho$ is
 half of the sum of all the positive roots, and  $2iH_{\rho}$ denotes
 the element in ${\frak t}$ corresponding to $2i\rho$ under the
 identification of ${\frak t}$ and ${\frak t}^*$  via the Killing
 form.  }
 \end{exam}

 \section{A cohomology pairing}
 \label{sec_poincare}

 Let $A$ be a Lie algebroid over $P$ with anchor map $\rho: A
 \rightarrow TP$. In this section, we show that  there is a natural
 pairing between $H^k(A)$, the $k$-th  Lie algebroid cohomology of $A$
 with trivial coefficients, and $H^{r-k}(A, Q_A)$, the $(r-k)$-th Lie
 algebroid cohomology  of $A$ with coefficients in $Q_A = \wedge^r A
 \ot \wedge^{top}T^*P$, where $r$ is the rank of $A$, and $0 \leq k
 \leq r$. For simplicity, we assume that $P$ is compact and
 orientable.  For a general $P$, we need to consider the ``compactly
 supported" Lie algebroid cohomology of $A$ and replace
 $\wedge^{top}T^*P$ by $|\wedge^{top}T^*P|$, the density bundle of $P$.

 Recall that
 \[
 H^k(A) ~ = ~ H^k(C, ~ d_{\ta}), \hspace{.4in} 
 H^k(A, ~ Q_A) ~ = ~ H^k(\tilde{C}, ~ \tilde{d}_{\ta}), 
 \]
 where $C = \oplus_{k=0}^{r}
 C^k = \oplus_{k=0}^{r} \Gamma(\wedge^k A^*)$ with $d_{\ta}$ given by
 (\ref{eq_dA}), and $\tilde{C} = C \ot \Gamma(Q_A)$ with 
 $ \tilde{d}_{\ta}: \tilde{C}^k \rightarrow \tilde{C}^{k+1}$
 given by
 \[
 \tilde{d}_{\ta} (\xi \ot s) ~ = ~ d_{\ta} \xi ~ \ot ~ s ~ + ~ (-1)^{|\xi|}
 \xi ~ \ot ~ Ds,
 \]
 where, for $s = X \ot \mu \in \Gamma(Q_A), ~ Ds \in \Gamma(A^*) \ot
 \Gamma(Q_A)$ maps a section $a$ of $A$ to the section
 \[
 D_a(s) ~ = ~ [a, ~ X] ~ \ot ~ \mu ~ + ~ X ~ \ot ~ L_{\rho(a)} \mu
 \]
 of $Q_A$ (see (\ref{eq_d-on-qa})). 

 Fix an orientation of $P$. Then integrating with respect to this
 orientation gives rise to the following  non-degenerate pairing
 between $C^k$ and $\tilde{C}^{r-k}$:
 \begin{equation}
 \label{eq_pairing}
 (\!( \xi, ~ \eta \ot X \ot \mu )\!) ~ {\stackrel{def}{=}} \int (\xi
 \wedge \eta, ~ X) \mu.
 \end{equation}

 \begin{thm}
 \label{thm_pairing}
 For $\xi \in C^{k-1}$ and $\eta \ot X \ot \mu \in  \tilde{C}^{r-k}$,
 \begin{equation}
 \label{eq_pairing-d}
 (\!( d_{\ta} \xi, ~ \eta \ot X \ot \mu)\!) ~ + ~  (-1)^{|\xi|}
 (\!(\xi, ~ \tilde{d}_{\ta} (\eta \ot X \ot \mu)  )\!) ~ = ~ 0.
 \end{equation}
 Consequently, there is an induced pairing between  $H^k (A)$ and
 $H^{n-k} (A, Q_A)$.
 \end{thm}

 \noindent
 {\bf Proof. }  Consider $c = (\xi \wedge \eta) \ot X \ot \mu \in
 \tilde{C}^{r-1}$.  We know that 
 \[
 \tilde{d}_{\ta} (c) ~ = ~  d_{\ta} \xi \ot (\eta \ot X \ot \mu) ~ + ~
 (-1)^{|\xi|} \xi \wedge \tilde{d}_{\ta} (\eta \ot X \ot \mu) \in
 \tilde{C}^r.
 \]
 Theorem  \ref{thm_pairing} now follows from the following Stokes'
 Theorem.  \qed

 \begin{thm}[Stokes' Theorem]
 \label{thm_stokes}
 Identify $\tilde{C}^r = \Gamma(\wedge^r A^* \ot \wedge^r A \ot
 \wedge^{top}T^*P)$ with the space of top-degree forms on $P$ by
 pairing the factors in $\wedge^r A^*$ and $\wedge^r A$ pointwise. We
 have, for  every $c = \xi \ot X \ot \mu \in \tilde{C}^{r-1}$,
 \begin{equation}
 \label{eq_a-a-r-1}
 \tilde{d}_{\ta}(c) ~ = ~ (-1)^{r-1} d( \rho(\xi \backl X) \backl \mu).
 \end{equation}
 Consequently,
 \begin{equation}
 \label{eq_int-0}
 \int_P \tilde{d}_{\ta}(c) ~ = ~ 0.
 \end{equation}
 \end{thm}

 \noindent
 {\bf Proof.} We only need to prove (\ref{eq_a-a-r-1}) locally.  Let
 $U$ be an open subset of $P$ over which both $\wedge^r A$ and
 $\wedge^{top} T^*P$ are trivial with nowhere vanishing  sections
 $X_0$ and $\mu_0$ respectively. Set $s_0 = X_0 \ot \mu_0$.  Let
 $\theta_0$ be the section of $A^*$ over $U$ such that 
 \[
  \theta_0 (a) s_0 ~ = ~ D_a s_0
 \]
 for every section $a$ of $A$. Write $X \ot \mu = f X_0 \ot \mu_0 =  f
 s_0$ over $U$. We have,
 \[
 \tilde{d}_{\ta} (c) ~ = ~ (d_{\ta} (f \xi) + \theta_0 \wedge f \xi, ~
 X_0) \mu_0.
 \]
 Set $a = f \xi \backl X_0$. It follows from the definitions of
 $\theta_0$ and the representation of $A$ on $Q_A$ that 
 \[ (\theta_0
 \wedge f \xi, ~ X_0 ) \mu_0 ~ = ~ - (d_{\ta} (f \xi), ~  X_0) \mu_0 ~
 +~ (-1)^{r-1} L_{\rho(a)} \mu_0.
 \]
 Hence,
 \[
 \tilde{d}_{\ta}(c) ~ = ~ (-1)^{r-1} L_{\rho(a)} \mu~ = ~  (-1)^{r-1}
 d (\rho(\xi \backl X) \backl \mu).
 \]
 This proves (\ref{eq_a-a-r-1}). It now follows from  Stokes' Theorem
 for de Rham cohomology that 
 \[
 \int_P \tilde{d}_{\ta}(c) ~ = ~ 0.
 \]
 \qed

 \begin{cor}
 \label{cor_pairing-on-poisson-homology}
 Let $P$ be an orientable Poisson manifold with a fixed orientation.
 For a differential form $\alpha$ and a compactly supported
 differential form $\beta$, define
 \[
 (\alpha, ~ \beta) ~ = ~ \int_P \alpha \wedge \beta.
 \]
 Let $\pdp: \Omega^k(P) \rightarrow \Omega^{k-1}(P)$ be the Koszul-Brylinski
 operator as defined by (\ref{eq_pp-1}).  Then, for any forms $\alpha$
 and $\beta$ with $|\alpha| + |\beta| - 1 = n = \dim P$, we have
 \begin{equation}
 \label{eq_pairing-on-poisson-homology}
 (\pdp \alpha, ~ \beta) ~ + ~ (-1)^{(|\alpha| - 1)}  (\alpha,~ \pdp \beta)
 ~ = ~ 0.
 \end{equation}
 Thus we get an induced pairing between the Poisson homology  spaces
 $H_k(P)$ and $H_{n-k}(P)$.
 \end{cor}

 \noindent
 {\bf Proof.} Fix a volume form $\mu_0$ of $P$. Let $\theta_0$ be the
 modular vector field defined by $\mu_0$.  Let $U \in \chi^{n -
 |\alpha|} (P)$ and  $V \in \chi^{n - |\beta|}(P)$ be such that
 \[
 \alpha ~ = ~ U \backl \mu_0, \hspace{.3in} \beta ~ = ~ V \backl \mu_0.
 \]
 let $(\!( ~ ~ )\!)$ be the pairing on $\chi(P)$ given by
 \[
 (\!( V_1, ~ V_2)\!) ~ = ~ \int_P (V_1 \wedge V_2, ~ \mu_0) \mu_0.
 \]
 Then we  know from Theorem \ref{thm_pairing} that 
 \[
 (\!([\pi, ~ U], ~ V) \!) ~ + ~ (-1)^{|U|}  (\!(U, ~ [\pi, V] + 2
 \theta_0 \wedge V)\!) ~ = ~ 0.
 \]
 It follows that 
 \begin{eqnarray*}
 & & (\!([\pi, ~ U] + \theta_0 \wedge U, ~ V) \!) ~ + ~ (-1)^{|U|}
 (\!(U, ~ [\pi, V] +  \theta_0 \wedge V)\!)\\ & = & (\!([\pi, ~ U], ~
 V) \!) ~ + ~ (-1)^{|U|}  (\!(U, ~ [\pi, V]) \!) ~ + ~ (\!( \theta_0
 \wedge U, ~ V)\!)  ~ + ~ (-1)^{|U|} (\!(U, ~  \theta_0 \wedge V)\!)
 \\ & = & (\!([\pi, ~ U], ~ V) \!) ~ + ~ (-1)^{|U|}  (\!(U, ~ [\pi,
 V]) \!)  ~ + ~ 2(-1)^{|U|} (\!(U, ~  \theta_0 \wedge V)\!) \\ & = &
 (\!([\pi, ~ U], ~ V) \!) ~ + ~ (-1)^{|U|}  (\!(U, ~ [\pi, V] + 2
 \theta_0 \wedge V)\!) \\ & = & 0.
 \end{eqnarray*}
 From Theorem \ref{thm_iso-bry}, we have 
 \beqa
 (\pdp \alpha, ~ \beta) & = & (-1)^{|U|+1}  (\!([\pi, ~ U] + \theta_0
 \wedge U, ~ V )\!)\\
 (\alpha, ~ \pdp \beta) & = & (-1)^{|V|+1} (\!( U, ~ [\pi, ~ V] +
 \theta_0 \wedge V) \!).
 \eeqa
 Thus we get (\ref{eq_pairing-on-poisson-homology}).  \qed
 
 We now turn to the discussion of the non-degeneracy of the  pairing
 $(\!( ~~, ~~)\!)$ between  $H^k(A)$ and $H^{r-k}(A, Q_A)$. We will
 see that it is not always non-degenerate. We first look at some
 familiar examples.

 \begin{exam}
 \label{exam_dual-g}
 {\em  When $A$ is a Lie algebra $\fg$ considered as a Lie algebroid
 over a one point space, we have $Q_A = \wedge^r \fg$, where $r = \dim
 \fg$.  The pairing $(\!( ~~, ~~)\!)$ in Theorem \ref{thm_pairing} is
 non-degenerate and it gives rise to an isomorphism
 \[
 H^k(\fg)^* ~ \cong ~ H^{r-k}(\fg, ~ \wedge^r \fg).
 \]
 This is the familiar Poincare duality for Lie algebra cohomology.  }
 \end{exam}

 \begin{exam}
 \label{exam_dual-tp}
 {\em Let $P$ be a compact orientable manifold and let $A = TP$ be 
 the tangent bundle Lie algebroid with  the
 identity anchor map. The line bundle $Q_A$ is trivial and so is the
 representation of $A$ on $Q_A$. The pairing $(\!( ~~, ~~)\!)$ in
 Theorem \ref{thm_pairing} is the  one obtained by integrating the
 wedge product of $k$-forms and $(n-k)$-forms, where $n = \dim
 P$. This, of course, is non-degenerate and gives the well-known
 Poincare duality for de Rham cohomology of $P$.  }
 \end{exam}

 \begin{exam}
 \label{exam_dual-0}
 {\em Let $A$ be an arbitrary vector bundle over $P$ of rank $r$.
 Consider $A$ as a Lie algebroid over $P$ with the zero Lie bracket on
 its sections and the zero anchor map to $TP$. Then $H^k(A)$ is the
 space of smooth sections of $\wedge^k A^*$, and $H^{r-k}(A, Q_A)$ is
 the space of smooth sections of 
 \[
 \wedge^{r-k}A^* \ot \wedge^r A \ot \wedge^{top}T^*P ~ \cong ~
 \wedge^k A \ot \wedge^{top}T^*P.
 \]
 The pairing between these two spaces, pairing the  elements in
 $\wedge^k A^*$ and in $\wedge^k A$ pointwise  and integrating the
 resulting top degree form on $P$ over $P$, is again non-degenerate in
 each argument.  }
 \end{exam}

 \begin{exam}
 \label{exam_dual-transformation}
 {\em We now consider the case when $A = P \times \fg$ is a
 transformation Lie algebroid defined by a Lie algebra  homomorphism
 $\rho: \fg \rightarrow \chi^1(P)$ from a Lie algebra $\fg$ to the Lie
 algebra of vector fields on $P$. See Example
 \ref{exam_transformation}. In this case,  both $H^{\bullet}(A)$ and
 $H^{\bullet}(A, Q_A)$ are natually isomorphic to certain Lie algebra
 cohomology spaces of $\fg$,  and the pairing $(\!( ~~, ~~)\!)$ in
 Theorem \ref{thm_pairing} is the one that occurs in the Poincare
 duality of Lie algebra cohomology (see, for example,
 \cite{kp:yellow}). More precisely, we have
 \begin{equation}
 \label{eq_cp-1}
 H^k(A) ~ \cong ~ H^k(\fg, ~ C^{\infty}(P)),
 \end{equation}
 where the right hand side is the Lie algebra cohomology of  $\fg$
 with coefficients in $C^{\infty}(P)$ considered as a $\fg$-module
 via $\rho$:
 \begin{equation}
 \label{eq_x-dot}
 x \cdot f ~ = ~ \rho(x) \cdot f, \hspace{.3in}  x \in \fg, ~ f \in
 C^{\infty}(P),
 \end{equation}
 and 
 \begin{equation}
 \label{eq_cp-2}
 H^{r-k}(A, Q_A) ~ \cong ~ H^{r-k}(\fg, ~ \wedge^r \fg \ot
 \Omega^{top}(P)),
 \end{equation}
 where $r = \dim \fg$ is the rank of $A$, and the right hand side is
 the Lie algebra cohomology of $\fg$ with coefficients  in the tensor
 product module $\wedge^r \fg \ot \Omega^{top}(P)$.  The space
 $\wedge^r \fg$ is equipped with the adjoint  representation of $\fg$,
 and $\Omega^{top}(P)$ is  equipped with the action of $\fg$ by Lie
 derivatives:
 \[
 x \cdot \mu ~ = ~ L_{\rho(x)} \mu.
 \]
 The isomorphisms (\ref{eq_cp-1}) and (\ref{eq_cp-2}) are easily seen
 to come from isomorphisms on the complex level: \beqa
 & & C^k ~ = ~ \Gamma(\wedge^k A^*) ~ \cong ~ \wedge^k \fgs \ot
 C^{\infty}(P)\\
 & & \tilde{C}^{r-k} ~ = ~ \Gamma(\wedge^{r-k} A^*  \ot \wedge^r A \ot
 \wedge^{top}T^*P) ~ \cong ~  \wedge^{r-k} \fgs  \ot \wedge^r \fg \ot
 \Omega^{top}(P).  \eeqa

 Consider the pairing between $C^{\infty}(P)$ and $\Omega^{top}(P)$
 given by
 \[
 (f, ~ \mu) ~ = ~ \int_P f \mu.
 \]
 For each $x \in \fg$, it follows from
 \[
 \int_P L_{\rho(x)}(f \mu) ~ = ~ 0
 \]
 that 
 \[
 (x \cdot f, ~ \mu) ~ + ~ (f, ~ x \cdot \mu) ~ = ~ 0.
 \]
 Thus, the two $\fg$-modules $C^{\infty}(P)$ and $\Omega^{top}(P)$ are
 contragradient to each other  with respect to the above
 pairing. Under the isomorphisms above, 
 the pairing $(\!( ~~ )\!)$ in Theorem
 \ref{thm_pairing} becomes the one between $\wedge^k \fgs \ot
 C^{\infty}(P)$ and $\wedge^{r-k} \fgs \ot \wedge^r \fg \ot
 \Omega^{top}(P)$ given by 
 \[
 (\xi \ot f, ~ \eta \ot X_0 \ot \mu) ~ = ~  (\xi \wedge \eta, ~ X_0)
 (f, ~ \mu).
 \]
 This is exactly the pairing that occurs  in
 the Poincare duality for Lie algebra cohomology.

 In this example, even though the pairing between the two
 $\fg$-modules $C^{\infty}(P)$ and $\Omega^{top}(P)$ is
 non-degenerate, they are not the full duals of each other,
 so we can not conclude that the induced  pairing on the cohomology
 spaces is non-degenerate. For example,
 \[
 H^0(\fg, C^{\infty}(P)) ~ = ~ C^{\infty}(P)^{\frak g}
 \]
 and
 \[
 H^r(\fg, ~ \wedge^r \fg \ot \Omega^{top}(P)) ~ \cong ~
 {\frac{\Omega^{top}(P)}{\fg \cdot \Omega^{top}(P)}}.
 \]
 These two spaces are not necessarily dual to each other, as we see
 from the following example: the vector field
 \[
 v_0 ~ = ~ (1 - e^{i \theta})^N {\frac{d}{d\theta}},
 \]
 where $N \geq 2$ is an integer,  defines an action of the
 $1$-dimensional Lie algebra $\fg = {\Bbb R}^1$ on $P = S^1$. The
 space $H^0(\fg, ~ C^{\infty}(P))$ is $1$-dimensional, while the
 space $H^r(\fg, ~ \wedge^r \fg \ot \Omega^{top}(P))$ is at least
 $N$-dimensional.  }
 \end{exam}

 \begin{exam}
 \label{exam_general-h0}
 {\em For a general Lie algebroid $A$ over $P$,  we have
 \[
 H^0(A) ~ \cong ~ \{f \in C^{\infty}(P): ~  \rho(a) \cdot f = 0, ~
 \forall a \in \Gamma(A)\}.
 \]
 and, from Stokes' Theorem,
 \[
 H^{r}(A, ~ Q_A) ~ \cong ~ {\frac{\Omega^{top}(P)}{\{ L_{\rho(a)} \mu,
 ~ a \ot \mu \in \Gamma(A) \ot \Omega^{top}(P)\} }}.
 \]
 If $P$ is not compact, we consider compactly supported  top-degree
 forms in $H^{r}(A, Q_A)$. The following is another  example where the
 two spaces $H^0$ and $H^r$ are not dual to each other: let $P = {\Bbb
 R}^2$ with the  Poisson structure given by
 \[
 \{x, ~ y\} ~ = ~ (x^2 + y^2)^N,
 \]
 where $N$ is an integer and $N \geq 3$. Let $A$ be $T^*P$ with the
 cotangent bundle Lie algebroid  defined by this Poisson
 structure. Again, $H^0(A)$ is $1$-dimensional and $H^r(A, Q_A)$ is at
 least $(N-1)$-dimensional.  }
 \end{exam}

 \section{The holomorphic case}
 \label{sec_holom}

 In this section, we extend our results to the holomorphic setting.
 Closely related results are studied in \cite{bz:outer} in different
 language, where questions of monodromy are also discussed.

 The theory of Lie algebroids can be developed in the setting of
 sheaves (see \cite{bb:bb} \cite{kt:sheaves}). Let $O_P$ be
 either the sheaf of  smooth functions on a manifold or the sheaf of
 holomorphic functions on a complex manifold (one can also consider
 the algebraic setting, but it will not be useful for us). Let $TP$ be
 the tangent sheaf  defined over $O_P$.  Note that $TP$ is a sheaf of
 Lie algebras over the scalars, where the Lie algebra structure on
 sections over an open set is the usual Lie bracket of vector fields. 

 A Lie algebroid $A$ over $O_P$ is a sheaf of $O_P$ modules together
 with 1) a Lie algebra structure on $A$ making $A$ into a sheaf of Lie
 algebras over the scalars, and 2) a homomorphism $\rho:A\to TP$ 
of sheaves of
 modules over $O_P$ and of Lie
 algebras, such that for $f \in O_P(U)$ and $\omega_{1},
 \omega_{2} \in A(U)$, the following derivation law holds:
 \[
 \{\omega_1, f \omega_2 \} = f \{\omega_1, \omega_2\} + (\rho
 (\omega_1) f) \omega_2.
 \]
  As before, $\rho$ is called the anchor map of the Lie algebroid. 
 For the rest of this section, 
 we will always assume that $A$ is locally free, so that it is
 the sheaf of sections of a vector bundle over $O_P.$ 

 A representation of $A$ is a sheaf $M$ of $A$ modules,  i.e., for
 every open set $U,$ the space of sections  $M(U)$ is a representation
 of $A(U)$ satisfying the properties given in Section  \ref{sec_intro}
 and the representations are compatible with restriction maps.  It
 follows as in Section  \ref{sec_modular} that the line bundle  $\qa =
 \wedge^{top}A \ot \wedge^{top}T^*P$ determines an $A$ module. 

 Consider the sheaf $\Omega^k(A)$ of sections of $\wedge^k A^*.$ The
 formula for $d_A$ given in Section \ref{sec_calculus} makes the
 sequence
 \[
 ... \to \Omega^{k-1}(A)\to \Omega^k(A)\to \Omega^{k+1}(A)\to ...
 \]
 into a complex of sheaves. We can take the hypercohomology of this
 complex of sheaves ${\hyperH} ^{\bullet}(\Omega^{\bullet}(A))$ (see
 \cite {gh:prin} p. $445$).  We define Lie algebroid cohomology
 \[
 H^{\bullet}(A,O_P)={\hyperH}^{\bullet}(\Omega^{\bullet}(A)).
 \]

 We also consider the cohomology sheaves ${\calH}^{\bullet} (A,O_P).$
 To define these sheaves, first consider the cohomology of the above
 complex on any open set. The assignment 
 \[
  U\mapsto H^i(\Omega^{\bullet}(A)(U),d_A)
 \]
 determines a presheaf. By definition,  ${\calH}^{\bullet}(A,O_P)$ is
 the sheaf associated to this presheaf.

 We will compute $H^{\bullet}(P,\calS)$ for a sheaf $\calS$ using Cech
 cohomology, so $H^i(A,O_P)$ is the total cohomology of the double
 complex $C^p(P,\Omega^q (A)).$ The usual filtrations give two first
 quadrant spectral sequences converging to $H^{\bullet} (A,O_P)$ \cite
 {gh:prin}.  The first has $E_2$ term $H^p(P,{\calH}^q(A, O_P))$ and
 the second has $E_2$ term $H^{q}(H^p(P,\Omega^\bullet (A))).$

 In the case where $O_P$ is smooth functions on $P,$ we have
 $H^p(P,\Omega^\bullet (A)))=0$ for $p>0$ since the sheaves
 $\Omega^{\bullet}(A)$ are soft.  It follows that in this case the
 second spectral sequence degenerates at the $E_2$ term, so that
 $H^q(A,O_P)=H^q(\Gamma(P, \Omega^{\bullet}(A))).$ Thus our
 hypercohomology definition of Lie algebroid cohomology agrees with
 the definition given in Section \ref{sec_calculus}. The same remarks
 apply when $O_P$ is the sheaf of holomorphic functions on a Stein
 manifold $P.$ 

 We return to the general case. We can compute Cech cohomology using
 an open cover $U_\alpha$ of $P$ such that all intersections of the 
 $U_\alpha$'s are contractible or empty. In particular, we may
 assume that on an intersection of open sets, every nowhere vanishing
 function has a logarithm. On such a cover, $Q_A$ is a trivial line
 bundle.  We can define a modular class $\Theta_A\in
 H^0(P,{\calH}^1(A,O_P))$ as follows. On each open set $U_\alpha$ in
 our cover, choose a nowhere vanishing section $s_{\alpha}$ of
 $Q_A|U_\alpha.$ Define a modular class on $U_\alpha$  by setting 
 \[
 D_as_{\alpha}=\Theta_{s_{\alpha}}(a)s_{\alpha}.
 \]
 As in Section \ref{sec_modular}, it follows that the class of
 $\Theta_{s_{\alpha}}\in {\calH}^1(A,O_P)$ is independent of the choice
 of $s_{\alpha}.$ We will denote it by $\Theta_\alpha .$ Moreover, 
 $\Theta_\alpha $ and $\Theta_\beta $
 agree on open sets $U_\alpha$ and $U_\beta$ in our cover. Indeed, on
 $U_\alpha \cap U_\beta,$ we know $s_{\alpha}=g_{\alpha\beta}s_{\beta}$
 for some nowhere vanishing function $g_{\alpha\beta}.$ Then it
 follows as above that the cohomology classes $\Theta_\alpha$ and
 $\Theta_\beta$ coincide. Thus we have a well-defined global section
 of the sheaf ${\calH}^1(A,O_P)$ which we denote by $\Theta_A.$ In
 fact, the functions $g_{\alpha \beta}$ are just the transition
 functions of $Q_A.$

 We wish to determine whether $\Theta_A$ determines a class in
 $H^1(A,O_P).$ We showed in the previous paragraph that $\Theta_A$
 determines a class in $H^0(P,{\calH}^1(A,O_P)),$ which is in
 $E_{2}^{0,1}$ for the first spectral sequence discussed
 above. $\Theta_A$ defines a class in $H^1(A,O_P)$ if it is
 annihilated by all differentials in the spectral sequence. All
 differentials $d_r$ for $r>2$ annihilate $\Theta_A$ for reasons of
 degree. Thus it suffices to determine if $d_2\Theta_A=0$ for the
 differential 
 \[
 d_2:H^0(P,{\calH}^1(A,O_P))\to H^2(P,{\calH}^0(A, O_P)).
 \]
  Denote the Cech differential by $\delta.$ By definition, for a class
 $r\in  H^0(P,{\calH}^1(A,O_P)),$ $d_2(r)=\delta c,$ where $c\in
 C^1(P,{\calH}^0(A,O_P))$ is chosen so $d_Ac=\delta r$ (see
 \cite{bt:forms}).

 \begin{prop}
 \label{prop_cclass}
 Let $i:{\boldC}_P\to {\calH}^0(A,O_P)$ be the inclusion of the
 constant sheaf into the sheaf of functions annihilated by $A.$ Then
 \[   
 d_2(\Theta_A)=i_*(2 \pi \sqrt{-1} c_1(Q_A)),
 \]
  where $c_1(Q_A)$ is the first Chern class of the line bundle $Q_A$
 in the holomorphic setting and 
 of the complexification of $Q_A$ in the real setting which is $0$.
 \end{prop}

 \noindent
 {\bf Proof.} We defined $\Theta_A$ by choosing a nowhere vanishing
 section $s_{\alpha}\in \Gamma (U_\alpha,Q_A).$ Then 
 \[ 
 \delta \Theta_A (a) (U_\alpha\cap U_\beta )= {{D_a(s_{\alpha })}\over
 {s_{\alpha }}}|_{U_\alpha \cap U_\beta} -  {{D_a(s_{\beta })}\over
 {s_{\beta }}}|_{U_\alpha \cap U_\beta} =d_{\scriptscriptstyle
 A}(log(g_{\alpha\beta}))(a)
 \]
 where the $g_{\alpha\beta}$ are the transition functions of $Q_A$
 discussed above. Then 
 \[ d_2 (\Theta_A)(U_{\alpha}\cap U_\beta\cap U_\gamma)=log(g_{\alpha\beta})+
 log(g_{\beta\gamma})-log(g_{\alpha\gamma})
 \]
 This is, up to a factor of $2 \pi \sqrt{-1},$ the usual 
 sheaf cohomology description of the first Chern
 class of the line bundle $Q_A.$ Since we regard this class as an
 element of $H^2(P,{\calH}^0(A,O_P)),$ we write it as $i_*(2 \pi \sqrt{-1}
 c_1(Q_A)).$
 \qed

 \begin{cor}
 \label{cor_cvan}
 Suppose $c_1(Q_A)=0.$ Then $d_2(\Theta_A)=0.$
 \end{cor}

 When $d_2(\Theta_A)=0,$ we denote the corresponding class in
 $H^1(A,O_P)$ by $\theta_A.$

 We remark that it follows from the proposition that  
 $d_2(\Theta_A)=0$ when $O_P$ is the sheaf of smooth
 functions. Indeed, in this case, the line bundle $Q_A$ has an
 underlying real structure so $c_1(Q_A)=0$ (\cite {bt:forms}). In the
 holomorphic setting, $Q_A$ is always trivial if $P$ is a Stein
 manifold. In addition, if $P$ admits a holomorphic symplectic
 structure, then $\wedge^{top}T^*P$ is trivial, so $c_1(Q_{T^*P})=0.$
 This is the case for $G/H,$ where $G$ is a complex reductive group
 and $H$ is a maximal torus, and a Poisson structure analogous to that
 of \cite {lu-we:poi} can be introduced.  In addition, in the case
 where the anchor map is zero, $d_2(\Theta_A)=0$ even though
 $c_1(Q_A)$ may not vanish. 

 \section{Appendix A: The adjoint ``representation'' of a Lie algebroid}
 \label{sec_adjoint}

 In this appendix, we will describe a construction which includes as
 special cases the adjoint representation of a Lie algebra and the
 (dual of the) flat ``Bott'' connection on the normal bundle to a
 foliation.  Given a Lie algebroid $A$, we will construct a
 ``representation up to homotopy'' of $A$ on the ``formal difference''
 $A \ominus TP$.  Taking the highest exterior power of this object
 will yield the representation of $A$ on $Q_A$ described in Section
 \ref{sec_modular}. That is all we will use of our construction in
 this paper, but we believe that it is interesting in its own right.

 The idea behind our construction is similar to that in the
 construction of the representation of $A$ on $Q_A$, namely that
 neither the Lie derivative action of $A$ on itself nor on $TP$ is a
 representation, but that ``the anomalies cancel.''  The following
 discussion makes this idea precise.

 \begin{dfn}
 \label{dfn_homotopy}
 {\em Let $A$ be a Lie algebroid over $P$ and let $(E,\del)$ be a
 bundle of $\zee_2$-graded complexes over $P$; i.e., $E$ is a
 $\zee_2$-graded vector bundle over $P$ and $\del$ is a bundle map of
 degree 1 with $\del ^2 = 0$.  A {\bf representation up to homotopy}
 of $A$ on $(E,\del)$ is an $\reals$-bilinear map
 \[
 \Gamma(A) ~ \times ~ \Gamma(E) \lrw \Gamma(E): ~~ a ~ \ot~  s \Map
 D_a s,
 \]
 such that the operators $D_a$ preserve the grading and commute with
 the action of $\del$ on sections, and such that the properties \beqa
 & & (2) ~~ D_a(fs) ~ = ~ f D_a s + (\rho(a)f) s;\\ & &  (3) ~~
 D_a(D_bs) - D_b(D_as) ~ = ~ D_{[a, b]} s.  \eeqa of a representation
 hold, while property \beqa & & (1) ~~ D_{fa}s ~ = ~ f D_a s \eeqa
 holds only up to homotopy, in the sense that for each $a \in
 \Gamma(A)$ and $f \in C^{\infty}(P)$ there is a bundle map
 $I(a,f):E\lrw E$ of degree 1 such that  \beqa & & (1') ~~ D_{fa}s ~
 = ~ f D_a s + I(a,f)\del s + \del I(a,f) s.  \eeqa }
 \end{dfn}

 \begin{rem}
 \label{rem_homotopy}
 {\em If we drop condition (3) from the definition of a
 representation, we get the definition of an $A$-{\bf connection} on a
 vector bundle $E$.  The difference of the left and right hand sides
 of (3) is then $C^{\infty}(P)$-linear in $a$, $b$, and $s$ and
 becomes the curvature tensor of the connection.  Similarly, we can
 drop condition (3) from Definition \ref{dfn_homotopy} to define an
 $A$-connection up to homotopy, and we could also weaken (3) to
 require flatness only up to homotopy in the definition of a
 representation.  Finally, one might put further conditions on 
 the trilinear expressions $I(a,f) s$ as one does in the the theory
 of strongly homotopy Lie algebras \cite{la-ma:strongly}. We have had
 no need to explore these options yet, though.  } 
 \end{rem}

 For a bundle $(E,\del)$ of complexes over $P$, the homology $H(E)$ is
 not a vector bundle unless $\del$ has constant rank.  We can still
 consider $H(E)$ as a $C^{\infty}(P)$-module by looking at the action
 of $\del$ on sections of $E$, though, so that the notion of a
 representation of $A$ on $H(E)$ still makes sense.  Also, we will
 define the {\bf determinant line bundle} $\wedge^{top}E$ of the
 graded vector bundle $E=E_0\oplus E_1$ to be $\wedge^{top}E_0^*
 \otimes \wedge^{top} E_1$.

 \begin{prop}
 \label{prop_homotopy}
 A representation up to homotopy of a Lie algebroid $A$ on a bundle
 $(E,\del)$ of $\zee_2$-graded complexes induces representations (not
 just up to homotopy) on the homology $H(E)$ and the  determinant line
 bundle $\wedge^{top}E$.  
 \end{prop}

 \noindent
 {\bf Proof.}  Since each $D_a$ commutes with $\del$, it induces an
 operator $H(D_a)$ on homology.  Properties (2) and (3) of the $D_a$'s
 are inherited by the induced operators.  To verify (1), we let $s$ be
 a cycle and find from $(1')$ that $D_{fa}s-f D_a s$ is a boundary.

 When $\del$ has constant rank, the determinant bundle of $H(E)$ is
 isomorphic to that of $E$, so the second part
 of the theorem follows from the first part in that special situation.
 For the general case, we will define and study the extended operators
 locally.  

 First of all, we note that a family of operators $D_a$ on sections of
 a vector bundle $V$ satisfying (2) and (3) can be extended in a
 unique way to a family of operators on all the tensors over $V$,
 satisfying the same identities, by requiring that the operators be
 derivations with respect to tensor product and commute with
 contractions.  This is done just is as usually done for Lie
 derivatives or covariant derivatives.  In particular, we obtain
 operators, also denoted by $D_a$, on $E^*$, $\wedge^{top}E_0^*$,
 $\wedge^{top} E_1$, and $\wedge^{top}E_0^* \otimes \wedge^{top} E_1$.

 For instance, the operators on $E^*$ are defined by 
 $$(D_a \omega,s)=\rho(a)\cdot(\omega,s) - (\omega,D_a s),$$ from
 which it follows that we have a representation up to homotopy on the
 dual complex $(E^*,\del^*)$, using the bundle maps $-I(a,f)^*$ to
 satisfy $(1')$.  

 To compute the operators on $\wedge^{top} E_1$, we use a local basis
 $t_1,\ldots,t_l$ for the sections of $E_1$.  By the derivation
 property,
 $$ D_a(t_1\wedge\ldots\wedge t_l)=\sum_{r=1}^l (t_1\wedge\ldots\wedge
 D_a t_r \wedge\ldots \wedge t_l).
 $$
 Then \beqa D_{fa}(t_1\wedge\ldots\wedge t_l) &=& \sum_{r=1}^l
 (t_1\wedge\ldots\wedge D_{fa} t_r \wedge\ldots \wedge t_l) \\ &=& f
 D_{a}(t_1\wedge\ldots\wedge t_l) \\ &+& \sum_{r=1}^l
 (t_1\wedge\ldots\wedge (I^{10}(a,f)\del^{01}+\del^{10} I^{01}(a,f))
 t_r \wedge\ldots \wedge t_l), \eeqa where the superscript $^{ij}$ on
 $I$ or $\del$ refers to the part of that bundle map which goes from
 $E_j$ to $E_i$.

 Expressing $(I^{10}(a,f)\del^{01}+\del^{10} I^{01}(a,f))t_r$ in terms
 of our basis, we see that all but its $r$'th component is annihilated
 by another factor in the wedge product, so that we get a simple
 expression in terms of a trace:
 \begin{equation}
 \label{eq_trace}
 D_{fa}(t_1\wedge\ldots\wedge t_l)= fD_{a}(t_1\wedge\ldots\wedge
 t_l)+\trace (I^{10}(a,f)\del^{01}+\del^{10}
 I^{01}(a,f))t_1\wedge\ldots\wedge t_l.
 \end{equation}
 Similarly, in terms of a basis $s_1^*,\ldots,s_k^*$ of local sections
 of $E_0^*$, we have:
 \begin{equation}
 \label{eq_tracestar}
 D_{fa}(s_1^* \wedge\ldots \wedge s_k^*)= fD_{a}(s_1^* \wedge\ldots
 s_k^*)-\trace  (I^{10*}(a,f)\del^{01*}+\del^{10*} I^{01*}(a,f))s_1^*
 \wedge\ldots s_k^*.
 \end{equation}
 Combining the two previous equations and using the derivation
 property with respect to tensor product, we obtain an expression for
 the behavior of the $D_a$'s operating on the determinant bundle:

 $$D_{fa}(s_1^* \wedge\ldots s_k^* \otimes t_1\wedge\ldots\wedge t_l)
    =fD_{a}(s_1^* \wedge\ldots s_k^* \otimes t_1\wedge\ldots\wedge
    t_l) +K(a,f)s_1^* \wedge\ldots s_k^* \otimes t_1\wedge\ldots\wedge
    t_l,
 $$
 where 
 $$K(a,f)=\trace (I^{10}(a,f)\del^{01}+\del^{10} I^{01}(a,f))- \trace
 	  (I^{10*}(a,f)\del^{01*}+\del^{10*} I^{01*}(a,f)).$$ Using
 	  the invariance of the trace under dualization and exchange
 	  of factors,   we can cancel all the terms in this expression
 	  to conclude that  $K(a,f)=0$, so that the operation of $A$
 	  on the determinant bundle satisfies $D_{fa}=fD_a$, so that
 	  we have an honest representation.   \qed

 With the general notion of representation up to homotopy at hand, we
 can consider our principal example.  Given the Lie algebroid $A$ with
 anchor $\rho:A\rightarrow TP$, we let $E_0=TP$, $E_1=A$,
 $\del^{01}=\rho$, and $\del^{10}=0$.  The homology of this little
 complex is the normal ``bundle'' $TP/\rho(A)$ to the orbits in degree
 0, and the isotropy ``bundle'' $\mbox{ker}~\rho$ in degree 1.  We call
 this the {\bf normal complex} of the Lie algebroid and denote it by
 $N(A)$.  The determinant bundle of the normal complex is precisely
 what we have called $Q_A$.

We define the operators $D_a$ on $\Gamma(E)$ for $a\in \Gamma(A)$ by 
$D_a b= [a,b]$ for $b \in \Gamma(A)$ and by
$D_a u= [\rho(a),u]$ for $u \in \Gamma(TP)$, where the last bracket is
the usual bracket of vector fields.  

\begin{prop}
\label{prop_normal}
The operators defined above form a representation up to homotopy of
 the Lie algebroid $A$ on $N(A)$.
 \end{prop}

\noindent
{\bf Proof.}
The fact that the $D_a$ are chain maps follows from the fact that
 $\rho$ defines a Lie algebroid homomorphism.  Properties (2) and (3)
 in the definition of a representation are standard facts in the
 differential calculus on Lie algebroids (see Section
 \ref{sec_calculus}).  For property $(1')$, we note as in Section
 \ref{sec_modular} that for the action on $A$ we have 
$$D_{fa} b = f D_a b - (\rho(b)\cdot f) a$$
 and for the action on $TP$ we have
 $$D_{fa}u = f D_a u - (u\cdot f) \rho(a).$$  
 If we define the homotopy operators by $I^{10}(a,f)(u)=-(df(u))a$ and
 $I^{01}(a,f)=0$, we see immediately from the two displayed equations
 above that $(1')$ is satisfied, so we have a representation up to homotopy.
\qed

When $A$ is a Lie algebra, the operators $D_a$ just give the
 adjoint representation, so we refer to them as the adjoint
``representation'' in the general case. When $\rho$ is the inclusion
of a subbundle of $TP$, the homology of the normal complex is just the
normal bundle to the corresponding foliation, and we recover the usual
flat connection along the leaves.  Finally, the associated
representation of $A$ on $Q_A$ may now be seen as the ``top
exterior power of the adjoint representation'' for a general Lie
algebroid.  

\section{Appendix B: The adjoint ``representation'' and modular class
 of a Lie groupoid}
\label{sec_groupoid}

Let $G$ be a Lie groupoid (i.e. a differentiable groupoid) over $P$.
We will denote its target and source maps by $\alpha$ and $\beta$, so
that the product $gh$ is defined whenever $\beta (g)=\alpha(h)$.
The groupoid analog of the Lie algebroid
``representation'' constructed in the previous section ought to be a
representation up to homotopy of $G$ on the
normal complex of its Lie algebroid $A$.  There are many possible
definitions of the notion of representation up to homotopy for a
groupoid, and we have not yet found an optimal one, so we will limit
ourselves here to a discussion of the particular case of the adjoint
``representation.''  

Recall that a representation of $G$
on a vector bundle
$\lambda:B\to P$ consists of a mapping $(g,r)\mapsto gr$ from 
$G\times_P B = \{(g,r)\in G\times B | \beta(g)=\lambda(r)\}$ to $B$
which is linear on fibres of $B$ and which satisfies the axioms:
\beqa 
& & (1) ~~ \lambda(gr)=\alpha(g);\\
& & (2) ~~ (gh)r = g(hr);\\
& & (3) ~~ er = r \mbox{ when $e$ is an identity element of $G$.}
\eeqa
In the case of Lie algebroids, the distinction between a
representation up to homotopy in our sense and an honest
representation is that $D_a s (x)$ can depend on the entire section
$a$ in the former case, while it depends only on $a(x)$ in the
latter.  In fact, for the adjoint ``representation'', $D_a s (x)$
depends only on the 1-jet of $a$ at $x$, and this carries over to the
groupoid case.

We will denote by $J^1 G$ the 1-jet prolongation groupoid of the Lie
groupoid $G$ over $P$.  The elements of $J^1 G$ are the 1-jets of
bi-sections of $G$, i.e.,  submanifolds of $G$ which project
diffeomorphically to $P$ under the source and target maps.  

Unlike the groupoid $G$ itself, $J^1 G$ has a natural
representation on the normal complex $N(A)$ of the Lie algebroid of $G$;
i.e., $J^1 G$ has representations on $A$ and $TP$ for which the anchor $\rho$
is an equivariant map.  One way to see this
is to consider elements of $J^1 G$ as special subsets of the
tangent bundle groupoid $TG$ over $TP$.  $TG$ acts on itself by left
translations, leaving invariant the part $T_P G$ of $TG$ lying over the
identity section of $G$.  Identifying $TP$ with a subbundle of $T_P G$
and $A$ with the normal bundle $T_P G /TP$ allows to obtain the
required representation of $J^1 G$.  

For example, if $G$ is the pair groupoid $P\times P$, an element
$(x,y)$ of $G$ does not naturally transport tangent vectors from $y$
to $x$, but an element of $J^1 G$ is precisely a vector space
isomorphism from $T_y P$ to $T_x P$.

We would like to make the action of $J^1 G$ descend to $G$ via the
natural ``forgetful'' projection $j:J^1 G \to G$.  As the example
above shows, this is not possible, but it turns out to be possible
``up to homotopy'' in the sense that, if $g'_1$ and $g'_2$ are two
elements of $J^1 G$ lying over the same $g\in G$, then the mappings
induced by the $g'_i$ between the complexes $N_{\alpha (g)}A$ and
$N_{\beta (g)}A$ are homotopic.  Proving the latter statement may be
reduced to the case where $g$ is an identity element at $x\in P$, in which
case the ``difference'' between $g'_1$ and $g'_2$ can be considered as
a linear map from $T_x P$ to the fibre $A_x$ of the Lie algebroid, and
this map (together with the zero map in the other direction) provides
the required homotopy operator.

As a consequence of this representation up to homotopy, the
representations of $J^1 G$ on $A$ and $TP$ descend to induce honest
representations of $G$ on the (generally singular) normal bundle
$TP/\rho(A)$ to the orbits and isotropy bundle $\mbox{ker}~\rho$ and
on the determinant line bundle $Q_A$.   In particular, the latter
representation defines an element of first cohomology of $G$ with
values in the multiplicative real numbers.  This is the modular class
of the groupoid.

In the situation of Example \ref{exam_foliation}, where $A$ is an integrable
subbundle of $TP$ we may take $G$ to be the holonomy groupoid of the
foliation, in which case we recover the linearized holonomy
representation of the foliation and the modular class of the foliation
in groupoid cohomology discussed in \cite{ya:modular} and Chapter IV
of \cite{mo-sc:global}.

 \end{document}